\documentclass[smallextended,referee,envcountsect,]{svjour3}
\smartqed

\usepackage{graphicx}
\usepackage{amsthm}
\usepackage{subfigure}
\usepackage{threeparttable}
\usepackage{amsmath}
\usepackage{color}
\usepackage{amssymb}
\usepackage{algorithmic}
\usepackage{algorithm}
\usepackage{array}
\usepackage{lipsum}
\usepackage{graphicx}
\usepackage{epstopdf}
\usepackage{ragged2e} 
\usepackage{booktabs,makecell, multirow, tabularx}
\usepackage[misc]{ifsym}
\newtheorem{assumption}{Assumption}
\numberwithin{example}{section} 
\journalname{JOTA}

\begin{document}

\title{A Smoothing Consensus-Based Optimization Algorithm for Nonsmooth Nonconvex Optimization}

\titlerunning{SCBO for nonsmooth nonconvex optimization}  

\author{Jiazhen Wei        \and
	Wei Bian}

\institute{Jiazhen Wei \at
		School of Mathematics, Harbin Institute of Technology\\
		 Harbin 150001, China \\
		jiazhenwei98@163.com
           \and
               Wei Bian,  Corresponding author  \at
              School of Mathematics, Harbin Institute of Technology\\
               Harbin 150001, China\\
              bianweilvse520@163.com
}

\date{Received: date / Accepted: date}

\maketitle

\begin{abstract}
Lately, a novel swarm intelligence model, namely the consensus-based optimization (CBO) algorithm, was introduced to deal with the global optimization problems. Limited by the conditions of Ito's formula, the convergence analysis of the previous CBO finite particle system mainly focuses on the problem with smooth objective function. With the help of smoothing method, this paper achieves a breakthrough by proposing an effective CBO algorithm for solving the global solution of a nonconvex, nonsmooth, and possible non-Lipschitz continuous minimization problem with theoretical analysis, which dose not rely on the mean-field limit. We indicate that the proposed algorithm exhibits a global consensus and converges to a common state with any initial data. Then, we give a more detailed error estimation on the objective function values along the state of the proposed algorithm towards the global minimum. Finally, some numerical examples are presented to illustrate the appreciable performance of the proposed method on solving the nonsmooth, nonconvex minimization problems.
\end{abstract}


\section{Introduction}
Global optimization has received increased attention across almost every field of economy, science and artificial intelligence.
For convex optimization problems, the fundamental result that a locally optimal solution is also a global
one is helpful for theoretical analysis and algorithm design. In addition, nonconvexity turns out to be ubiquitous in large-scale
optimization problems \cite{DemoA2021,HertleinAn2019,RuszA2021}. Common examples include big data analysis, signal processing and resource allocation. Although nonconvexity
enables us to accurately describe practical problems, it usually brings great challenges to its research on theory and algorithms. Finding a global solution is not easy if there are more than one local solution, where the main difficulty is how to escape from the local minimizers and choose the right one. Moreover, it has been shown that a range of nonconvex optimization problems are not only NP-hard to solve optimally, but NP-hard to solve approximately as well \cite{MekaRank2008}.

Besides the common used deterministic method, stochastic (or heuristic) method is another class of currently being adopted basic methods for global optimization. A stochastic method is to 
generate sampling points by random perturbation and use them for
nonlinear local search, where the classical stochastic optimization methods often include Genetic Algorithms (GAs) \cite{Grefenstette1986Optimization},
Tabu Search \cite{SiaryTabu2000}, Simulated Annealing (SA) \cite{Hwang1988Simulated}, Swarm Intelligence (SI) \cite{Tarasewich2002Swarm,YangSwarm2018}, etc.
By simulating the behavior of animal groups in nature, SI algorithm uses the information of exchange and cooperation among groups to achieve the global optimization by simple and limited interactions among individuals. Comparing with the traditional optimization
algorithm, it can search for the global optimal solution of some complex optimization problems more quickly and effectively.

In this paper, we consider the following optimization problem:
\begin{equation}\label{problem}
	\min_{x\in\mathbb{R}^d} f(x),
\end{equation}
where $f:\mathbb{R}^d\rightarrow\mathbb{R}_{+}$ is a continuous function that admits a unique global minimizer $x^*\in\mathbb{R}^d$.
Let $x^i(t)=(x_1^i(t),\ldots,x_d^i(t))^\mathrm{T}\in\mathbb{R}^d$ be the coordinate process of the $i$-th sample point at time $t$. Pinnau et al. \cite{PinnauA2016}
proposed a novel first-order stochastic SI algorithm called consensus-based optimization (CBO) algorithm, which can be expressed by the following stochastic differential equations $(i=1,\ldots,N)$:
\begin{equation}\label{cboo}
	\begin{cases}
		\displaystyle \mathrm{d}x^i(t)=-\lambda(x^i(t)-\upsilon_f)H^\varepsilon \left(f(x^i(t))-f(\upsilon_f)\right)\mathrm{d}t+\sqrt{2}\sigma\|x^i(t)-\upsilon_f\|\mathrm{d}W^i(t),\\
		\displaystyle \upsilon_f(t)=\dfrac{\sum_{i=1}^Nx^i(t)\mathrm{e}^{-\beta f(x^i(t))}}{\sum_{i=1}^N\mathrm{e}^{-\beta f(x^i(t))}},\\
	\end{cases}
\end{equation}
where $\lambda>0$, $\sigma\geq 0$, $\beta$ is a positive constant corresponding to the reciproal of temperature in statistical physics,
$H^\varepsilon$ denotes a smooth regularization of the Heaviside function and $\{W^i(t)\}_{i=1}^N$ are independent Brownian motions. They studied
its convergence behavior from a mean-field point of view. Then, a simple version of the CBO algorithm in \cite{PinnauA2016} was proposed in \cite{CarrilloAn2016}, where the system is as follows:
\begin{equation}\label{cbo2}
	\begin{cases}
		\displaystyle \mathrm{d}x^i(t)=-\lambda(x^i(t)-\upsilon_f)\mathrm{d}t+\sigma\|x^i(t)-\upsilon_f\|\mathrm{d}W^i(t), ~i=1,\ldots,N,\\
		\displaystyle \upsilon_f(t)=\dfrac{\sum_{i=1}^Nx^i(t)\mathrm{e}^{-\beta f(x^i(t))}}{\sum_{i=1}^N\mathrm{e}^{-\beta f(x^i(t))}}.\\
	\end{cases}
\end{equation}
Under some appropriate conditions, the authors in \cite{CarrilloAn2016} showed that the algorithm in \eqref{cbo2} converges to the global
minimizer of \eqref{problem} in the mean-field sense. Subsequently, in order to eliminate the dependence of drift rate on dimensions
and make CBO algorithm more competitive in high-dimensional optimization problems, the diffusion term in \eqref{cbo2} was replaced by a component-wise diffusion in \cite{CarrilloA2019}, where the algorithm is modeled by
\begin{equation}\label{cbo}
	\begin{cases}
		\displaystyle \mathrm{d}x^i(t)=-\lambda(x^i(t)-\upsilon_f)\mathrm{d}t+\sigma\sum\limits_{l=1}^d(x_l^i(t)-\upsilon_f)\mathrm{d}W^i_l(t)e_l,i=1,\ldots,N,\\
		\displaystyle \upsilon_f=\dfrac{\sum_{i=1}^Nx^i(t)\mathrm{e}^{-\beta f(x^i(t))}}{\sum_{i=1}^N\mathrm{e}^{-\beta f(x^i(t))}},\\
	\end{cases}
\end{equation}
where $e_l$ is the column vector in $\mathbb{R}^d$ with $1$ at the $l$-th element and $0$ for the others. In both continuous and semi-discrete cases, the authors proved that the mean-field limit of the algorithm in \eqref{cbo} converges to a good approximation of the global minimizer, and the parameters are independent of the dimensionality. In addition, the authors in \cite{CarrilloA2019} utilized the random mini-batch strategy to calculate the
weighted mean, which reduces the computational cost and makes it be more suitable for high dimensional optimization problems. By drawing on \cite{CarrilloA2019}, Ha et al. \cite{HaConvergence2020} gave a direct proof to the convergence of a simplified model of \eqref{cbo} in the particle level. Here, the last summand term of \eqref{cbo} was replaced by $$\sigma\sum\limits_{l=1}^d(x_l^i(t)-\upsilon_f)\mathrm{d}W_l(t)e_l,$$ which means that it supposes that all particles are influenced by a shared environmental noise.
Recently, due to the promising theoretical results and numerical performance of the CBO algorithm established in the global optimization, more attention has been paid to the study of some variants of the CBO algorithm \cite{BorghiConstrained2023,Chenchene2023,FornasierConsenseus2021,HaConvergence2021,Huang2023,2020Consensus}.

Although the continuous CBO algorithm has been widely studied, except for \cite{HaConvergence2020}, its properties were still investigated on corresponding mean-field limit, a Fokker-Planck equation, in the other literatures.
\cite{FornasierConsensus2022} devised a novel wholistic proof to the convergence of CBO algorithm \eqref{cboo} to the global minimizer of a class of locally Lipschitz continuous functions in mean-field law. \cite{FornasierConsensus2024} proposed a CBO algorithm variant by introducing an additional truncation in the noise term. They also proved the convergence in expectation under some mild assumptions on the objective function and the initialization. In \cite{JoseFedCBO2023}, the authors proposed a new FedCBO algorithm that can be used to solve a class of non-convex optimization problems arising in practical clustered federated learning settings. They also gave the convergence analysis in the mean-field regime.
	Considering that the Fokker-Planck equation is derived by letting $N \rightarrow \infty$, the convergence of the Fokker-Planck equation to the global minimizer in mean-field law can shed lights on the properties of CBO model, but it can not rigorously imply the convergence of CBO algorithm per se. So it was stated in \cite{JoseFedCBO2023} that directly studying the finite particle system without passing to the mean-field limit is of particular interest. Additionally, the convergence analysis of CBO algorithm in \cite{HaConvergence2020} needs the objective function $f$ to be $\mathcal{C}_b^2(\mathbb{R}^d)$, which means that $f$ is a $\mathcal{C}^2$ function on $\mathbb{R}^d$ with bounded derivatives up to second-order. Moreover, note that the estimate in \cite[Theorem 4.1]{HaConvergence2020} is an error estimate for the CBO algorithm with a fixed parameter $\beta$, which can not directly imply that the error can be small enough. To sum up, motivated by the above research, we study the convergence of the related CBO algorithm without resorting to the corresponding mean-field model, and the global optimization problem, 
	whose objective function $f$ can be relaxed drastically to a nonsmooth or even non-Lipschitz continuous function in this paper.

 In Table \ref{tablesummary}, we summarize and compare the convergence results of our proposed algorithm with other aforementioned CBO algorithms from the following two aspects:
\begin{itemize}
	\item the analyzed model: finite particle model or mean-field limit model;
	\item assumptions on the objective function $f$.
\end{itemize}
\begin{table}[htbp]
\caption{Summary of the convergence results of our algorithm and other related CBO methods.} \label{tablesummary}
	\begin{center}
		\renewcommand\arraystretch{1.4}
		\resizebox{\textwidth}{!}{
			\begin{threeparttable}
				\begin{tabular}{|c|c|c|} \hline
					Algorithm&Analyzed model&Assumptions on $f$ \\
					\hline
					\multirow{3}{*}{SCBO [this paper]}&\multirow{3}{*}{finite particle}&continuous;\\
					&&\makecell{its smoothing function \\can be constructed}\\
					\hline
					\multirow{2}{*}{\makecell{\eqref{cbo2} with common \\noise \cite{HaConvergence2020}}}&\multirow{2}{*}{finite particle}&$\mathcal{C}^2$;\\
					&&$\nabla^2 f$ is bounded\\
					\hline
					\multirow{2}{*}{CBO in \eqref{cboo} \cite{FornasierConsensus2022}}&\multirow{2}{*}{mean-field limit}&local Lipschitz continuous;\\
					&&inverse continuity condition\\
					\hline
					\multirow{3}{*}{\makecell{CBO with truncated \\noise \cite{FornasierConsensus2024}}}&\multirow{3}{*}{mean-field limit}&local Lipschitz continuous;\\
					&&inverse continuity condition;\\
					&&the growth of $f$ is bounded\\
					\hline
					\multirow{4}{*}{FedCBO \cite{JoseFedCBO2023}}&\multirow{4}{*}{mean-field limit}&local Lipschitz continuous;\\
					&&the growth of $f$ is bounded;\\
					&&$\nabla f$ is Lipschitz continuous;\\
					&&$\nabla f$ is bounded\\
					\hline
				\end{tabular}
		\end{threeparttable}}
\end{center}
\end{table}

Though the CBO algorithm is independent of the gradient of $f$, its nonsmoothness makes us unable to apply the essential Ito's formula to its error estimate.
Therefore, we introduce the smoothing method, a well-known method for solving nonsmooth optimization problems for many decades
\cite{BianA2020,ChenSmoothing2012,ChenPenalty2016}, to overcome the nonsmoothness of $f$ in this paper. Taking advantage of the smoothing method, we can solve a nonsmooth optimization problem by a sequence of optimization problems with smooth objective functions. 
Therefore, in this paper, we only require that the objective function $f$ in \eqref{problem} is continuous. Actually, for a continuous function, there is a framework to construct its smooth approximation functions by convolution \cite{ChenSmoothing2012,HiriartConvex1993,RockafellarVariational1998}. Here, we would like to emphasize that $f$ in \eqref{problem} can be neither convex nor smooth, or even not locally Lipschitz continuous necessarily, which can be applied to a broader class of applications \cite{BianNeural2014,ChenA2018,HuangRobust2018}.
We construct a variant of the CBO model in \eqref{cbo}, in which the original function $f$ is replaced by the smoothing function of $f$ with a positive smoothing parameter.
First, we prove that the new algorithm exhibits a global consensus for any initial data and has a common consensus state almost surely. Then,
thanks to the designed update rules of the smoothing parameter, we provide a sufficient condition on the system parameters and initial data to guarantee the function values at the consensus state being a good
approximation of the global minimum. Finally, some numerical experiments are given to verify our theoretical results and show the effectiveness of the proposed algorithm in finding the global solution of the nonsmooth nonconvex optimization problems. In conclusion, the research in this paper is a complement and improvement to the theoretical analysis of the CBO algorithm while ensuring good performance on numerical experiments.

The structure of this paper is outlined as follows.
In Section \ref{sec:2}, we introduce some definitions and preliminary results which will be needed in this paper.
In Section \ref{sec:3}, we propose a variant of the CBO algorithm to solve \eqref{problem}. Then we study the emergence of global consensus of
the solution process to the proposed algorithm.
In Section \ref{sec:4}, we study the emergence of a common consensus state and the error analysis of the proposed algorithm to minimization problem \eqref{problem}.
In Section \ref{sec:5}, some numerical examples are given to support the theoretical results of this paper.
In Section \ref{sec:6}, we draw the conclusions.

\textbf{Notation:}~Denote $\mathbb{R}^d$ the $d$-dimensional real-valued vector space and $\mathbb{R}^d_+=[0,+\infty)^d$. Given column vectors $x=(x_1,\ldots,x_d)^\mathrm{T}$
and $y=(y_1,\ldots,y_d)^\mathrm{T}$, $x\cdot y=\sum_{i=1}^d x_iy_i$ is the scalar product of $x$ and $y$. 
$\mathcal{C}^2(\mathbb{R}^d)$ indicates the set of twice continuously differentiable functions on $\mathbb{R}^d$.
For random variables $X$ and $Y$ in the probability space $(\Omega,\mathcal{F},\mathbb{P})$, $X\sim Y$ means both $X$ and $Y$ follow the same distribution.
$\mathbb{E}[X]$ indicates the expectation of $X$, and $\mathbb{P}\{A\}$ indicates the probability of event $A\in\mathcal{F}$. $\|\cdot\|$ stands for the $\ell^2$-norm of a vector and $\|\cdot\|_2$ denotes the spectral norm of a matrix.
$e_l$ is the column vector in $\mathbb{R}^d$ with $1$ at the $l$-th element and $0$ for the others.
Parameter $\delta_{kl}$ is defined by $\delta_{kl}=1$ if $k=l$ and $\delta_{kl}=0$ if $k\neq l$. $\mathcal{P}_{ac}(\mathbb{R}^d)$ denotes the space of Borel probability measures
that are absolutely continuous w.r.t the Lebesgue measure on $\mathbb{R}^d$. 

\section{Preliminaries}
\label{sec:2}
In this section, we introduce some necessary definitions and preliminary results. We refer the readers to \cite{BerntStochastic1985,RossStochastic1995} for more details,
in which some mathematical preliminaries of stochastic process, Ito integral, martingale, stochastic differential equations are introduced.

First, we recall some definitions in stochastic process, including the   global consensus in $L^p$ and global consensus almost surely.
\begin{definition}[\cite{HaConvergence2021}]
	Let $(\Omega,\mathcal{F},\mathbb{P})$ be a probability space, and $X_t:=\{x^i(t)=(x_1^i(t),\ldots,x_d^i(t))^\mathrm{T}\}_{i=1}^N$ be a stochastic process  on it.
	\begin{itemize}
		\item [{\rm (i)}] We call the stochastic process $X_t$ exhibits a global consensus in $L^p$ with $p\geq1$, if the following $L^p$-convergence holds:
		$$\lim\limits_{t\rightarrow\infty}\mathbb{E}\left[ \|x^i(t)-x^j(t)\|^p\right] =0,~~\forall i,j=1,\ldots,N.$$
		\item [{\rm (ii)}] We call the stochastic process $X_t$ exhibits a global consensus almost surely (a.s.), if for almost all $\omega\in\Omega$ and $i,j=1,\ldots,N$, the sample path $x^i(t,\omega)-x^j(t,\omega)$ tends to zero asymptotically, i.e.
		$$\mathbb{P}\left\{\lim\limits_{t\rightarrow\infty}\|x^i(t)-x^j(t)\|=0\right\}=1,~~\forall i,j=1,\ldots,N.$$
	\end{itemize}
\end{definition}
Next, some preliminaries on stochastic process are presented.
Note that if $X_t$ is a martingale, then $X_t$ has a right-continuous modification which is a submartingale. Therefore, for a continuous martingale, there is a corresponding  martingale convergence theorem as follows.
\begin{proposition}[\cite{2013Brownian}]\label{marcon}
	Let $X_t$ be a supermartingale with right-continuous sample paths. Assume that the collection $\{X_t\}_{t\geq 0}$ is bounded in $L^1$. Then there
	exists a random variable $X_\infty\in L^1$ such that
	$$\lim_{t\rightarrow\infty}X_t=X_\infty,~~{\rm a.s.}.$$
\end{proposition}
\begin{proposition}[\cite{GrigoriosStochastic2014}]\label{ito}
	For the Ito stochastic integral $\int_0^tf(s)\mathrm{d}W(s)$, it has
	\begin{itemize}
		\item [{\rm (i)}] $\mathbb{E}\left[ \int_0^tf(s)\mathrm{d}W(s)\right] =0;$
		\item [{\rm (ii)}] (Ito isometry) $\mathbb{E}\left[\left(\int_0^tf(s)\mathrm{d}W(s)\right)^2\right]=\mathbb{E}\left[\int_0^tf^2(s)\mathrm{d}s\right].$
	\end{itemize}
\end{proposition}
\begin{proposition}[\cite{CrowLognormal1988}]\label{log-normal}
	A Lognormal distribution is a continuous distribution of random variable $X$ whose logarithm is normally distributed. And if the random variable $Y:=\ln X$ follows the normal distribution with expectation $\alpha$ and variance $\beta^2$, that is $Y\sim \mathcal{N}(\alpha,\beta^2)$, then 
	$\mathbb{E}[X]=\mathrm{e}^{\left(\alpha+\frac{\beta^2}{2}\right)}$.
\end{proposition}

Thanks to the above proposition, we can get the expectation of $\mathrm{e}^Y$, if random variable $Y\sim \mathcal{N}(\alpha, \beta^2)$. 

Ito's formula is one of the most important techniques for the theoretical analysis of this paper, which is introduced in what follows.
\begin{proposition}[Ito's formula \cite{BerntStochastic1985}]
	Assume that $x(t)=(x_1(t),\dots,x_d(t))^\mathrm{T}$ is a $d$-dimensional Ito's process satisfying
	$$\mathrm{d}x(t)=u\mathrm{d}t+v\mathrm{d}W(t),$$
	where
	\begin{equation*}
			u=\left(\begin{array}{c}
				u_1\\
				\vdots\\
				u_d\\
			\end{array}\right),~
			v=\left(\begin{array}{ccc}
				v_{11} & \cdots & v_{1m}\\
				\vdots & & \vdots\\
				v_{d1} & \cdots & v_{dm}\\
			\end{array}\right),~
			\mathrm{d}W(t)=\left(\begin{array}{c}
				\mathrm{d}W_1(t)\\
				\vdots\\
				\mathrm{d}W_m(t)\\
			\end{array}\right).
		\end{equation*}
		Let $g(x,t)=(g_1(x,t),\ldots,g_p(x,t))^\mathrm{T}$ be a $\mathcal{C}^2(\mathbb{R}^d \times [0,\infty))$ function. Then the process $Y(\omega,t)=g(x(t),t)$ is again an Ito's process, whose $q$-th component $Y_q$, is given by
		$$dY_q=\frac{\partial g_q(x,t)}{\partial t}\mathrm{d}t+\sum\limits_{l=1}^{d}\frac{\partial g_q(x,t)}{\partial x_l}\mathrm{d}x_l(t)+\frac{1}{2}\sum\limits_{l,k=1}^{d}\frac{\partial^2 g_q(x,t)}{\partial x_l \partial x_k}\mathrm{d}x_l(t)\cdot \mathrm{d}x_k(t)$$
		with
		\begin{equation}\label{dd}
			\mathrm{d}W_l(t)\cdot \mathrm{d}W_k(t)=\delta_{lk}\mathrm{d}t \quad \mbox{and} \quad
			\mathrm{d}t\cdot \mathrm{d}W_l(t)=\mathrm{d}W_l(t)\cdot \mathrm{d}t=\mathrm{d}t\cdot \mathrm{d}t=0.
		\end{equation}
	\end{proposition}
	
	By Laplace principle \cite{DemboLarge1998}, we further have the following result.
	\begin{proposition}[\cite{PinnauA2016}]\label{Laplace}
		Assume that $f:\mathbb{R}^d\rightarrow\mathbb{R}_+$ is bounded and attains the global minimum at the unique point $x^*$, then for any probability measure $\rho\in\mathcal{P}_{ac}(\mathbb{R}^d)$, it holds
		\begin{equation}\label{laplace}
			\lim\limits_{\beta\rightarrow\infty}\left(-\frac{1}{\beta}\log\left(\int_{\mathbb{R}^d}\mathrm{e}^{-\beta f(x)}\mathrm{d}\rho(x)\right)\right)=f(x^*).
		\end{equation}
	\end{proposition}
	\begin{remark}
		By \eqref{laplace}, it is obvious to know
		$$\lim_{\beta\rightarrow\infty}-\frac{1}{\beta}\log\mathbb{E}\left[ \mathrm{e}^{-\beta f(x(\omega))}\right] =f(x^*).$$
	\end{remark}
	
	Many smooth approximations for nonsmooth functions have been developed over several decades. Smoothing methods make it possible to overcome the nonsmoothness of the objective function, so that the Ito's formula can be applied to the corresponding smoothing functions. In this paper, the smoothing function to the objective function in \eqref{problem} is defined as follows.
	\begin{definition}[\cite{BianWorst2013}]\label{def-sf}
		For the objective function $f$ in \eqref{problem},
		we call $\tilde{f}:\mathbb{R}^d\times(0,\bar{\mu}]\rightarrow\mathbb{R}$ with $\bar{\mu}>0$ a smoothing function of it,
		if the following conditions hold:
		\begin{itemize}
			\item [{\rm (i)}] for any fixed $\mu\in(0,\bar{\mu}]$, $\tilde{f}(\cdot,\mu)$ is twice continuously differentiable in
			$\mathbb{R}^d$, and for any fixed $x\in\mathbb{R}^d$, $\tilde{f}(x,\cdot)$ is differentiable on $(0,\bar{\mu}]$;
			\item [{\rm (ii)}] for any fixed $x\in\mathbb{R}^d$, $\lim_{z\rightarrow x,\mu\downarrow0}\tilde{f}(z,\mu  )=f(x)$;
			\item [{\rm (iii)}] for any bounded set $\mathcal{X}\subseteq\mathbb{R}^d$, there exist positive constants $\kappa$ and $\eta$, and $q\in[0,1)$ such that for any $x\in\mathcal{X}$ and $\mu\in(0,\bar{\mu}]$, it holds
			\begin{equation}\label{eq-mu1}
				\left|\nabla_{\mu}\tilde{f}(x,\mu)\right|\leq\kappa\mu^{-q} \quad \mbox{and} \quad \left\|\nabla^2_{xx}\tilde{f}(x,\mu)\right\|_2\leq \eta\mu^{-q-1}.
			\end{equation}
		\end{itemize}
	\end{definition}
	Definition \ref{def-sf} (ii) and (iii) imply that
	\begin{equation}\label{eq-mu3}
		\left|\tilde{f}(x,\mu)-f(x)\right|\leq\kappa\mu^{1-q},~~\forall x\in\mathcal{X},\, \mu\in(0,\bar{\mu}].
	\end{equation}
	\begin{remark}\label{q}
		For a continuous function $f:\mathbb{R}^d\rightarrow\mathbb{R}$, we can construct a smoothing function of it based on convolution \cite{HiriartConvex1993,RockafellarVariational1998}. Consider a sequence of bounded, measurable functions $\psi^\mu\geq 0$ with $\int_{\mathbb{R}^d}\psi^\mu(z)\mathrm{d}z=1$ such that the sets $B^\mu=\{z|\psi^\mu(z)>0\}$ form a bounded sequence converging to $\{0\}$ as $\mu$ converging to $0$. Define function $\tilde{f}:\mathbb{R}^d\times\mathbb{R}_+\rightarrow\mathbb{R}$ by \begin{equation}\label{smoothing1}
			\tilde{f}(x,\mu):=\int_{\mathbb{R}^d}f(x-z)\psi^\mu(z)\mathrm{d}z=-\int_{\mathbb{R}^d}f(z)\psi^\mu(x-z)\mathrm{d}z.
		\end{equation}
	  From \eqref{smoothing1}, we can easily find that if $\psi^\mu$ is $r$th continuously differentiable, then $\tilde{f}(\cdot,\mu)$ is $(r+1)$th continuously differentiable. 
	  A common used method for constructing $\psi^{\mu}$ is to start from a compact set $B$ of the origin and a continuous function $\psi:B\rightarrow[0,\infty]$ with $\int_{B}\psi(z)\mathrm{d}z=1$, and to define $\psi^{\mu}$ to be $\mu^{-d}\psi(z/\mu)$ when $z\in \mu B$, and $0$ when  $z\not\in \mu B$ \cite[Example 7.19]{RockafellarVariational1998}. 
	  With this definition on $\psi^{\mu}$, \eqref{smoothing1} can be expressed by
	  \begin{equation}\label{smoothing22}
			\tilde{f}(x,\mu)=-\int_{\mathbb{R}^d}f(z)\mu^{-d}\psi((x-z)/\mu)\mathrm{d}z=\int_{\mathbb{R}^d}f(x-\mu{z})\psi(z)\mathrm{d}z,
		\end{equation}
		which satisfies Definition \ref{def-sf}-(i) if $\psi$ is continuously differentiable on $\mathbb{R}^d$, and satisfies Definition \ref{def-sf}-(ii) by \cite[Theorem 9.67]{RockafellarVariational1998}. Based on \eqref{smoothing22}, when $f$ is locally Lipschtz continuous, we can have the following estimations.
		\begin{itemize} 
		\item If $\int_{\mathbb{R}^d}\|z\|\psi(z)\mathrm{d}z\leq\bar{\rho}<+\infty$, then
		$$\nabla_{\mu}\tilde{f}(x,\mu)=\int_{\mathbb{R}^d\setminus{D}_x}\nabla f(x-\mu{z})^{\mathrm{T}}z\psi(z)\mathrm{d}z,$$
		where $D_x$ is a set of measure $0$ such that $f(x-\mu z)$ is differentiable with respect to $z$ for any $z\not\in D_x$. 
		Thus, $|\nabla_{\mu}\tilde{f}(x,\mu)|\leq\kappa$, $\forall x\in\mathcal{X}$, $\mu\in(0,\bar{\mu}]$ with $\kappa=L\bar{\rho}$, where $L$ is a Lipschitz constant of $f$ over set $\mathcal{X}+\bar{\mu}B$.
		\item If $\psi$ is second order continuously differentiable and $\int_{\mathbb{R}^d}\|\nabla\psi(z)\|\mathrm{d}z\leq\bar{\eta}$, then
		$$\nabla^2_{xx}\tilde{f}(x,\mu)=-\int_{\mathbb{R}^d}f(z)\mu^{-d-2}\nabla^2\psi((x-z)/\mu)\mathrm{d}z.$$
		Thus, $\|\nabla^2_{xx}\tilde{f}(x,\mu)\|\leq\eta\mu^{-1}$ with $\eta=L\bar{\eta}$.
		\end{itemize} 
		Thus, if $f$ is locally Lipschitz continuous on $\mathbb{R}^d$, we can always construct a smoothing function of it based on \eqref{smoothing22} with $q=0$. Moreover, when $h:\mathbb{R}^d\rightarrow\mathbb{R}$ is second order continuously differentiable and we can give a smoothing function of $g$, denoted by $\tilde{g}$, then $h+\tilde{g}$ is a smoothing function of $h+g$. 
				
		Moreover, we use a widely used example to illustrate this exactly. The plus function $t_+$ has been used in a variety of important areas throughout science and engineering, including penalty functions, complementarity problems, signal reconstruction, optimal control, stochastic equilibrium, spherical approximations, and so forth. Many nonsmooth functions can be reformulated by using the plus function $t_+$, such as $|t|$, $\max(t,y)$, $\min(t,y)$ and so on. Let $\rho: \mathbb{R}\rightarrow \mathbb{R}_+$ be a continuous density function satisfying $\rho(s)=\rho(-s)$ and $\kappa:=\int_{-\infty}^\infty|s|\rho(s)\mathrm{d}s<\infty$.
		Let
		\begin{equation*}
			\phi(t,\mu):=\int_{-\infty}^\infty(t-\mu s)_+\rho(s)\mathrm{d}s,
		\end{equation*}
		which satisfies 
		$$\nabla_{\mu}\phi(t,\mu)=\int^{+\infty}_{t/\mu}s\rho(s)\mathrm{d}s, \,\nabla_t\phi(t,\mu)=\int_{-\infty}^{t/\mu}\rho(s)\mathrm{d}s,\, \nabla_t^2\phi(t,\mu)=\rho(t/\mu)/\mu.$$ Thus, 
		$\phi(t,\mu)$ is a smoothing function of $t_+$ with $q=0$.  
		Based on $\phi$, we can construct smoothing functions for a class of nonsmooth functions which can be expressed by composition of the plus function $t_+$ \cite{ChenSmoothing2012}. And in \cite[Lemma 1]{ChenSmoothing2012}, it has been shown that, for any $p\in(0,1]$, $\psi(t,\mu)^p=\left( \phi(t,\mu)+\phi(-t,\mu)\right)^p$ is the smoothing function of $|t|^p$ with $q=1-p$. We refer to \cite{BianNeural2014,ChenSmoothing2012} for more details on this topic. Moreover, if $f:=|g(x)|^p$ with second order continuously differentiable function $g:\mathbb{R}^d\rightarrow \mathbb{R}$ and $p\in(0,1]$, then $\tilde{f}(x,\mu)=\psi(g(x),\mu)^p$ is a smoothing function of $f$ with $q=1-p$. And we can notice that the smoothing function defined in Definition \ref{def-sf} has the additivity, that is if $\tilde{f}_1$ and $\tilde{f}_2$ are smoothing functions of $f_1$ and $f_2$, respectively, then $\tilde{f}_1+\tilde{f}_2$ is a smoothing function of $f_1+f_2$.
		\end{remark}
		\begin{example}
			Consider the following function
			\begin{equation}\label{sfproblem}
				f(x)=\sum_{i=1}^m|x_i|^p
			\end{equation}
			with $0<p\leq 1$ and $x_i\in\mathbb{R}, i=1,\ldots,m$. Choose
			\begin{equation}\label{smoothing}
				\varphi_1(s,\mu)=\mu\ln\left( 2+\mathrm{e}^{-\frac{s}{\mu}}+\mathrm{e}^{\frac{s}{\mu}}\right) 
			\end{equation}
			as a smoothing function of $|s|$. Then
			$$\tilde{f}(x,\mu)=\sum_{i=1}^m\varphi_1(x_i,\mu)^p$$
			is a smoothing function of $f$ with $\kappa=mp(\ln4)^p$, $\eta=2mp(\ln4)^{p-1}+mp(1-p)(\ln4)^{p-2}$ and $q=1-p$.
		\end{example}
			\begin{example}
				We can also choose
				\begin{equation}\label{smoothing2}
					\varphi_2(s,\mu)=\sqrt{s^2+4\mu^2}
				\end{equation}
				as a smoothing function of $|s|$. Then 
				$$\tilde{f}(x,\mu)=\sum_{i=1}^m\varphi_2(x_i,\mu)^p$$
				is a smoothing function of $f$ in \eqref{sfproblem} with $\kappa=4mp2^{p-2}$, $\eta=8mp\left(\frac{8}{2-p}+4 \right)^{\frac{p-4}{2}}+mp2^{p-2}  $ and $q=1-p$.
		\end{example}
	
\section{Global Consensus}
\label{sec:3}

In this section, we first propose a variant of CBO algorithm with a smoothing function of $f$ in Subsection \ref{sec:3.1}.
Then, we provide some sufficient conditions on the system parameters and initial data to let our algorithm lead to a global consensus in Subsection \ref{sec:3.2}.
\subsection{A variant of CBO Algorithm with Smoothing Function}\label{sec:3.1}
To solve problem \eqref{problem}, now based on Definition \ref{def-sf},
we set $\mu:[0,\infty)\rightarrow(0,\bar{\mu}]$ as a continuously differentiable and monotone decreasing function satisfying $\lim_{t\rightarrow\infty}\mu(t)=0$.
Inspired by the CBO algorithm in \eqref{cbo}, we propose the Smoothing Consensus-Based Optimization (SCBO) algorithm described in Algorithm \ref{alg}.
\begin{algorithm}
	\normalsize
	\caption{SCBO algorithm}
	\label{alg}
	\begin{algorithmic}
		\STATE{\begin{align}
				&\mathrm{d}x^i(t)=-\lambda(x^i(t)-\bar{x}^*(t))\mathrm{d}t+\sigma\sum\limits_{l=1}^d(x^i_l(t)-\bar{x}^*_l(t))\mathrm{d}W_l(t)e_l,~i=1,\ldots,N.\label{alg-1} \\ 
				&\bar{x}^*(t)=(\bar{x}^*_1(t),\ldots,\bar{x}^*_d(t))^\mathrm{T}:=\dfrac{\sum_{i=1}^Nx^i(t)\mathrm{e}^{-\beta \tilde{f}(x^i(t),\mu(t))}}{\sum_{i=1}^N\mathrm{e}^{-\beta \tilde{f}(x^i(t),\mu(t))}}.\label{alg-2} \\ 
				&x^i(t)\Big|_{t=0}=x^i(0). \nonumber
			\end{align}
		}
	\end{algorithmic}
\end{algorithm}

In \eqref{alg-1}, drift rate $\lambda>0$, noise intensity $\sigma\geq 0$, and $\beta$ is a given positive constant. The one-dimensional Brownian motions $\{W_l(t)\}_{l=1}^d$ are assumed to be i.i.d and satisfy the following relations:
$$\mathbb{E}[W_l(t)]=0,~~l=1,\ldots,d,\quad \mathbb{E}[W_{l}(t)W_{k}(t)]=\delta_{kl}t,~~1\leq l,k\leq d.$$
Here we emphasize that the external noise is the same for all particles in the SCBO algorithm, which is the same as in \cite{HaConvergence2020}.

\subsection{Global Consensus}\label{sec:3.2}
Global consensus means that the particles will adjust their behaviors to reach a common state through the mutual cooperations and communications. From the view point of opinion dynamics, the emergence of global consensus also shows that all particles will share a common opinion. In this subsection, we will prove that for any initial data, the SCBO algorithm in Algorithm \ref{alg} exhibits a global consensus almost surely, which is important for the subsequent convergence analysis.

Considering the continuous process for $x^i$ and $x^j$ in \eqref{alg-1}, we have
\begin{equation}\label{cboij}
	\begin{cases}
		\displaystyle \mathrm{d}x^i(t)=-\lambda(x^i(t)-\bar{x}^*(t))\mathrm{d}t+\sigma\sum\limits_{l=1}^d(x^i_l(t)-\bar{x}^*_l(t))\mathrm{d}W_l(t)e_l,&t>0,\\
		\displaystyle \mathrm{d}x^j(t)=-\lambda(x^j(t)-\bar{x}^*(t))\mathrm{d}t+\sigma\sum\limits_{l=1}^d(x^j_l(t)-\bar{x}^*_l(t))\mathrm{d}W_l(t)e_l,&t>0.\\
	\end{cases}
\end{equation}

Subtracting the two equations in \eqref{cboij} to see that $x^i(t)-x^j(t)$ satisfies
\begin{equation}\label{dij}
	\mathrm{d}(x^i(t)-x^j(t))=-\lambda(x^i(t)-x^j(t))\mathrm{d}t+\sigma\sum\limits_{l=1}^d(x^i_l(t)-x^j_l(t))\mathrm{d}W_l(t)e_l.
\end{equation}
It follows from the $l$-th component of \eqref{dij} that $x^{ij}_l(t):=x^i_l(t)-x^j_l(t)$ satisfies
\begin{equation}\label{dijl}
	\begin{cases}
		\displaystyle \mathrm{d}x^{ij}_l(t)=-\lambda x^{ij}_l(t)\mathrm{d}t+\sigma x^{ij}_l(t)\mathrm{d}W_l(t),~~t>0,\\
		\displaystyle x^{ij}_l(t)\Big|_{t=0}=x^i_l(0)-x^j_l(0).\\
	\end{cases}
\end{equation}
Next, the following lemma gives the quadratic variation of $x^i_k(t)$ and $x^i_l(t)$.
\begin{lemma}\label{lem4.3}
	Let $\{x^i(t)\}_{i=1}^N$ be a solution process to the SCBO algorithm. Then the quadratic variation of $x^i_k(t)$ and $x^i_l(t)$ is given as follows:
	\begin{equation*}
		\mathrm{d}x^i_k(t)\cdot \mathrm{d}x^i_l(t)=
		\left\{\begin{array}{ll}
			\sigma^2(x^i_k(t)-\bar{x}^*_k(t))^2\mathrm{d}t,& \mbox{if}~~ k=l,\\
			0,& \mbox{if}~~ k\neq l.
		\end{array} \right.
	\end{equation*}
\end{lemma}
{\noindent\it Proof} It follows from \eqref{alg-1} that the $l$-th component of $x^i(t)$ satisfies
	\begin{equation}\label{dil}
		\mathrm{d}x^i_l(t)=-\lambda(x^i_l(t)-\bar{x}^*_l(t))\mathrm{d}t+\sigma(x^i_l(t)-\bar{x}^*_l(t))\mathrm{d}W_l(t).
	\end{equation}
	Then we use \eqref{dd} and \eqref{dil} to get
	\begin{equation*}
		\begin{aligned}
			\mathrm{d}x^i_k(t)\cdot \mathrm{d}x^i_l(t)
			=&\sigma^2\delta_{kl}(x^i_k(t)-\bar{x}^*_k(t))(x^i_l(t)-\bar{x}^*_l(t))\mathrm{d}t\\
			=&\left\{\begin{array}{ll}
				\sigma^2(x^i_k(t)-\bar{x}^*_k(t))^2\mathrm{d}t,& \mbox{if}~~ k=l,\\
				0,& \mbox{if}~~ k\neq l.
			\end{array} \right.
		\end{aligned}
	\end{equation*}\qed

Now, applying Ito's formula for $\ln x^{ij}_l(t)$, and using \eqref{dijl} gives
\begin{equation}\label{dlnijl}
	\resizebox{0.9\hsize}{!}{$\begin{aligned}
			&\mathrm{d}\ln x^{ij}_l(t)
			=\frac{\mathrm{d}x^{ij}_l(t)}{x^{ij}_l(t)}-\frac{1}{2(x^{ij}_l(t))^2}\mathrm{d}x^{ij}_l(t)\cdot \mathrm{d}x^{ij}_l(t)\\
			=&-\lambda \mathrm{d}t+\sigma \mathrm{d}W_l(t)\\
			&-\frac{1}{2(x^{ij}_l(t))^2}\left(-\lambda x^{ij}_l(t)\mathrm{d}t+\sigma x^{ij}_l(t)\mathrm{d}W_l(t)\right)
			\cdot\left(-\lambda x^{ij}_l(t)\mathrm{d}t+\sigma x^{ij}_l(t)\mathrm{d}W_l(t)\right)\\
			=&-\left(\lambda+\frac{\sigma^2}{2}\right)\mathrm{d}t+\sigma \mathrm{d}W_l(t),
		\end{aligned}$}
\end{equation}
where the second equality follows from \eqref{dijl} and the last equality uses the quadratic variation relations in \eqref{dd}.

Integrating the relation in \eqref{dlnijl} from $0$ to $t$ gives
\begin{equation}\label{ijl}
	x^{ij}_l(t)=x^{ij}_l(0)\exp\left(-\left(\lambda+\frac{\sigma^2}{2}\right)t+\sigma W_l(t)\right),~~t\geq 0,~l=1,\ldots,d.
\end{equation}
Then, the explicit formula of \eqref{ijl} implies the following conclusions.
\begin{theorem}\label{consensus1}
	Let $\{x^i(t)\}_{i=1}^N$ be a solution process of the SCBO algorithm. Then, it holds that
	\begin{itemize}
		\item [{\rm(i)}] the global consensus occurs a.s., i.e. for $i,j=1,\ldots,N$ and $l=1,\ldots,d$, it holds
		$$\lim\limits_{t\rightarrow\infty}| x^i_l(t)-x^j_l(t)|=0,~~{\rm a.s.};$$
		\item [{\rm(ii)}] for $i,j=1,\ldots,N$, it holds
		$$\mathbb{E}\left[ \|x^i(t)-x^j(t)\|^2\right] =\mathrm{e}^{-(2\lambda-\sigma^2)t}\mathbb{E}\left[ \|x^i(0)-x^j(0)\|^2\right] ,~~t\geq0.$$
		Hence, if $2\lambda>\sigma^2$, the global consensus in $L^2$ emerges a.s..
	\end{itemize} 
\end{theorem}
{\noindent\it Proof}
	(i) Recall the law of iterated logarithm of the Brownian motion \cite{AchimProbability2008}:
	\begin{equation}\label{logarithm}
		\limsup\limits_{t\rightarrow\infty}\frac{| W_l(t)|}{\sqrt{2t\log\log t}}=1,~~{\rm a.s.}.
	\end{equation}
	From \eqref{ijl} and the definition of $x^{ij}_l(t)$, we have
	\begin{equation}\label{xij}
		x^i_l(t)-x^j_l(t)=(x^i_l(0)-x^j_l(0))\exp\left(-\left(\lambda+\frac{\sigma^2}{2}\right)t+\sigma W_l(t)\right),~~\forall t\geq0.
	\end{equation}
	And by \eqref{logarithm}, note that for $t\rightarrow\infty$, the linear negative term $-(\lambda+\frac{\sigma^2}{2})t$ in $t$ is certainly dominant compared to the Brownian term $\sigma W_l(t)$.
	Thus, we obtain
	$$\lim_{t\rightarrow\infty}| x^i_l(t)-x^j_l(t)|=0,~~{\rm a.s.}.$$
	
	(ii) We square the both sides of \eqref{xij} to find
	\begin{equation}\label{square}
		(x^i_l(t)-x^j_l(t))^2=(x^i_l(0)-x^j_l(0))^2\exp\left(-(2\lambda+\sigma^2)t+2\sigma W_l(t)\right),~~\forall t\geq0.
	\end{equation}
	Then, summing up \eqref{square} for $l$ from $1$ to $d$, one has
	\begin{equation}\label{sum}
		\|x^i(t)-x^j(t)\|^2=\mathrm{e}^{-(2\lambda+\sigma^2)t}\sum_{l=1}^d\mathrm{e}^{2\sigma W_l(t)}(x^i_l(0)-x^j_l(0))^2,~~\forall t\geq0.
	\end{equation}
	By Proposition \ref{log-normal} and the fact that $W_l(t)$ follows the Gaussian distribution with mean $0$ and variance $t$, we obtain
	$$2\sigma W_l(t)\sim \mathcal{N}(0,4\sigma^2t)\quad\Rightarrow \quad\mathbb{E}\left[ \mathrm{e}^{2\sigma W_l(t)}\right] =\mathrm{e}^{2\sigma^2t}.$$
	Then, taking expectation on the both sides of \eqref{sum} gives
	\begin{equation}\label{l2}
		\begin{aligned}
			\mathbb{E}\left[ \|x^i(t)-x^j(t)\|^2\right] =&\mathrm{e}^{-(2\lambda+\sigma^2)t}\sum_{l=1}^d\mathbb{E}\left[ \mathrm{e}^{2\sigma W_l(t)}(x^i_l(0)-x^j_l(0))^2\right] \\
			=&\mathrm{e}^{-(2\lambda+\sigma^2)t}\sum_{l=1}^d\mathrm{e}^{2\sigma^2 t}\mathbb{E}\left[ (x^i_l(0)-x^j_l(0))^2\right] \\
			=&\mathrm{e}^{-(2\lambda-\sigma^2)t}\mathbb{E}\left[ \|x^i(0)-x^j(0)\|^2\right] ,~~\forall t\geq0.
		\end{aligned}
	\end{equation}
	Finally, if $2\lambda>\sigma^2$, we let $t$ tend to $\infty$  in \eqref{l2} to get
	$$\lim\limits_{t\rightarrow\infty}\mathbb{E}\left[ \|x^i(t)-x^j(t)\|^2\right] =0,~~\forall i,j=1,\ldots,N,$$
	which means the global consensus in $L^2$ emerges for the SCBO algorithm.
\qed

\section{Convergence Analysis}
\label{sec:4}

In this section, we provide the convergence analysis of the SCBO algorithm.
\subsection{Common Consensus State} \label{sec:4.1}

We have shown the emergence of global consensus for the SCBO algorithm in Section \ref{sec:3.2}.
But this does not mean that the solution tends to a common consensus state.
In this subsection, we will discuss the existence of the common consensus state for the SCBO algorithm.

First, we introduce a point based on the state of the SCBO algorithm, called the ensemble average and defined by:
$$\bar{x}(t)=(\bar{x}_1(t),\ldots,\bar{x}_d(t))^\mathrm{T}:=\frac{1}{N}\sum\limits_{i=1}^Nx^i(t).$$

The following elementary lemma, which will be crucially used later, gives some relationships between $x^i(t)$, $\bar{x}(t)$ and $\bar{x}^*(t)$.
\begin{lemma}\label{3}
	Let $\{x^i(t)\}_{i=1}^N$ be a solution process of the SCBO algorithm. Then, the following estimates hold.
	\begin{itemize}
		\item [{\rm (i)}]$\|x^i(t)-\bar{x}(t)\|^2=\sum\limits_{l=1}^d(x^i_l(0)-\bar{x}_l(0))^2\exp\left(-(2\lambda+\sigma^2)t+2\sigma W_l(t)\right);$
		\item [{\rm (ii)}]$\|\bar{x}(t)-\bar{x}^*(t)\|^2\leq \max_{1\leq i\leq N}\|x^i(t)-\bar{x}(t)\|^2;$
		\item [{\rm (iii)}]$\dfrac{1}{N}\sum\limits_{i=1}^N\mathbb{E}\left[ \|x^i(t)-\bar{x}^*(t)\|^2\right] \leq 2\mathrm{e}^{-(2\lambda-\sigma^2)t}\sum\limits_{l=1}^d\mathbb{E}\left[\max_{1\leq i\leq N}(x^i_l(0)-\bar{x}_l(0))^2\right].$
	\end{itemize}
\end{lemma}
{\noindent\it Proof}
	(i) Following from \eqref{alg-1} and the definition of $\bar{x}(t)$, we derive that
	\begin{equation}\label{bar}
		\mathrm{d}\bar{x}(t)=-\lambda(\bar{x}(t)-\bar{x}^*(t))\mathrm{d}t+\sigma\sum\limits_{l=1}^d(\bar{x}_l(t)-\bar{x}^*_l(t))\mathrm{d}W_l(t)e_l.
	\end{equation}
	Subtracting \eqref{alg-1} from \eqref{bar} shows
	\begin{equation*}
		\mathrm{d}(x^i(t)-\bar{x}(t))=-\lambda(x^i(t)-\bar{x}(t))\mathrm{d}t+\sigma\sum\limits_{l=1}^d(x^i_l(t)-\bar{x}_l(t))\mathrm{d}W_l(t)e_l,
	\end{equation*}
	whose $l$-th component is
	\begin{equation*}
		\mathrm{d}(x^i_l(t)-\bar{x}_l(t))=-\lambda(x^i_l(t)-\bar{x}_l(t))\mathrm{d}t+\sigma(x^i_l(t)-\bar{x}_l(t))\mathrm{d}W_l(t).
	\end{equation*}
	Combining with the analogous process of \eqref{dij}-\eqref{ijl}, one has
	\begin{equation}\label{e1}
		x^i_l(t)-\bar{x}_l(t)=(x^i_l(0)-\bar{x}_l(0))\exp\left(-\left(\lambda+\frac{\sigma^2}{2}\right)t+\sigma W_l(t)\right).
	\end{equation}
	Squaring the both sides of \eqref{e1} and summing up it over $l=1,\ldots,d$ yields that
	$$\|x^i(t)-\bar{x}(t)\|^2=\sum_{l=1}^d(x^i_l(0)-\bar{x}_l(0))^2\exp\left(-(2\lambda+\sigma^2)t+2\sigma W_l(t)\right).$$
	
	(ii) By the definition of $\bar{x}^*(t)$ in \eqref{alg-2}, it is obvious that
	$$\|\bar{x}(t)-\bar{x}^*(t)\|^2=\left\|\bar{x}(t)-\dfrac{\sum_{i=1}^Nx^i(t)\mathrm{e}^{-\beta \tilde{f}(x^i(t),\mu(t))}}{\sum_{i=1}^N\mathrm{e}^{-\beta \tilde{f}(x^i(t),\mu(t))}}\right\|^2\leq \max\limits_{1\leq i\leq N}\|\bar{x}(t)-x^i(t)\|^2.$$
	
	(iii) Note that
	\begin{equation}\label{iii}
		\begin{aligned}
			&\dfrac{1}{N}\sum\limits_{i=1}^N\|x^i(t)-\bar{x}^*(t)\|^2\\
			=&\dfrac{1}{N}\sum\limits_{i=1}^N\|x^i(t)-\bar{x}(t)+\bar{x}(t)-\bar{x}^*(t)\|^2\\
			=&\dfrac{1}{N}\sum\limits_{i=1}^N\left(\|x^i(t)-\bar{x}(t)\|^2+\|\bar{x}(t)-\bar{x}^*(t)\|^2\right)\\
			\leq& 2\max\limits_{1\leq i\leq N}\|x^i(t)-\bar{x}(t)\|^2\\
			=&2\max\limits_{1\leq i\leq N} \left(\sum\limits_{l=1}^d(x^i_l(0)-\bar{x}_l(0))^2 \exp\left(-(2\lambda+\sigma^2)t+2\sigma W_l(t)\right) \right)\\
			\leq&2\sum\limits_{l=1}^d\left(\max\limits_{1\leq i\leq N}(x^i_l(0)-\bar{x}_l(0))^2\right)\exp\left(-(2\lambda+\sigma^2)t+2\sigma W_l(t)\right),\\
		\end{aligned}
	\end{equation}
	where the second equality uses the fact that $\frac{1}{N}\sum_{i=1}^{N}(x^i(t)-\bar{x}(t))=0$, the first inequality follows from (ii) and the last equality follows from (i).
	Thus, by taking expectation on the both sides of \eqref{iii}, we get
	\begin{equation}\label{i}
		\begin{aligned}
			&\dfrac{1}{N}\sum\limits_{i=1}^N\mathbb{E}\left[ \|x^i(t)-\bar{x}^*(t)\|^2\right] \\
			\leq&
			2\sum\limits_{l=1}^d\left(\mathbb{E} \left[ \max\limits_{1\leq i\leq N}(x^i_l(0)-\bar{x}_l(0))^2\right]\right) \mathbb{E}\left[ \exp\left(-(2\lambda+\sigma^2)t+2\sigma W_l(t)\right)\right] \\
			=&2\mathrm{e}^{-(2\lambda-\sigma^2)t}\sum_{l=1}^d\left( \mathbb{E}\left[ \max\limits_{1\leq i\leq N}(x^i_l(0)-\bar{x}_l(0))^2\right]\right) ,\\
		\end{aligned}
	\end{equation}
	where the fact that $\mathbb{E}\left[ \mathrm{e}^{2\sigma W_l(t)}\right] =\mathrm{e}^{2\sigma^2t}$ is used in the last equality.
\qed
\begin{theorem}\label{consensus}
	Let $\{x^i(t)\}_{i=1}^N$ be a solution process of the SCBO algorithm. Suppose the system parameters $\lambda$ and $\sigma$ in it satisfy $2\lambda>\sigma^2$,
	then there exists a common random vector $x_\infty$ such that
	$$\lim\limits_{t\rightarrow\infty}x^i(t)=x_\infty, ~~{\rm a.s.},\quad \forall 1\leq i\leq N.$$
\end{theorem}

{\noindent\it Proof}
	Integrating \eqref{alg-1} from $0$ to $t$ shows
	\begin{equation}\label{I12}
		x^i_l(t)=x^i_l(0)-\lambda\mathcal{I}_l^i(t)+\sigma\tilde{\mathcal{I}}_l^i(t), ~i=1,\ldots,N,~l=1,\ldots,d,
	\end{equation}
	where
	$$\mathcal{I}_l^i(t)=\int_0^t(x^i_l(s)-\bar{x}^*_l(s))\mathrm{d}s \quad \mbox{and} \quad  \tilde{\mathcal{I}}_l^i(t)=\int_0^t(x^i_l(s)-\bar{x}^*_l(s))\mathrm{d}W_l(s).$$
	Next, we show the almost surely convergence of the terms $\mathcal{I}_l^i(t)$ and $\tilde{\mathcal{I}}_l^i(t)$, $i=1,\ldots,N$, $l=1,\ldots,d$, separately.
	
	On the one hand, owing to \eqref{iii}, we have
	\begin{equation*}
			\begin{aligned}
				&\left|x^i_l(t)-\bar{x}^*_l(t)\right|\\
				\leq&\sqrt{\sum\limits_{i=1}^N\|x^i(t)-\bar{x}^*(t)\|^2}\\
				\leq&\sqrt{2N\sum\limits_{l=1}^d\left(\max\limits_{1\leq i\leq N}(x^i_l(0)-\bar{x}_l(0))^2\right)\exp\left(-(2\lambda+\sigma^2)t+2\sigma W_l(t)\right)}.\\
			\end{aligned}
	\end{equation*}
	Combining the above inequality with \eqref{logarithm}, we obtain that there exist positive random variables $C_i:=C_i(\omega),i=1,2$ such that
	\begin{equation}\label{C12}
		| x^i_l(t)-\bar{x}^*_l(t)|\leq C_1\mathrm{e}^{-C_2t},\quad{\rm a.s.}~~\omega\in\Omega,\quad\forall t\geq 0.
	\end{equation}
	Then, we set
	$$\mathcal{J}_l^i(t):=\mathcal{I}_l^i(t)-\int_0^tC_1\mathrm{e}^{-C_2s}\mathrm{d}s=\int_0^t\left(x^i_l(s)-\bar{x}^*_l(s)-C_1\mathrm{e}^{-C_2s}\right)\mathrm{d}s.$$
	As a result of \eqref{C12}, the integrand above is nonpositive a.s., and $\mathcal{J}_l^i$ is consequently non-increasing in $t$ a.s..
	In addition, note that
	$$|\mathcal{J}_l^i(t)|\leq \int_0^t2C_1\mathrm{e}^{-C_2s}\mathrm{d}s=\frac{2C_1}{C_2}\left(1-\mathrm{e}^{-C_2t}\right)\leq\frac{2C_1}{C_2},$$
	where the first inequality follows from \eqref{C12}.
	Since $\mathcal{J}_l^i$ is monotone decreasing and bounded from below along sample paths, we obtain
	$$\exists\mathcal{J}_\infty^{i,l}:=\lim\limits_{t\rightarrow\infty}\mathcal{J}_l^i(t)=\lim\limits_{t\rightarrow\infty}\left(\mathcal{I}_l^i(t)-\int_0^tC_1\mathrm{e}^{-C_2s}\mathrm{d}s\right),~~{\rm a.s.}.$$
	This yields
	$$\lim\limits_{t\rightarrow\infty}\mathcal{I}_l^i(t)=\mathcal{J}_\infty^{i,l}+\frac{C_1}{C_2},~~{\rm a.s.}.$$
	
	On the other hand, we obtain:
	\begin{equation}\label{mar}
		\begin{aligned}
			&\mathbb{E}\left[ \int_0^t(x^i_l(s)-\bar{x}^*_l(s))^2\mathrm{d}s\right] \\
			\leq&\int_0^t\sum\limits_{i=1}^N\mathbb{E}\left[ \|x^i(s)-\bar{x}^*(s)\|^2\right] \mathrm{d}s\\
			\leq&2N\left(\int_0^t\mathrm{e}^{-(2\lambda-\sigma^2)s}\mathrm{d}s\right)\sum_{l=1}^d\left(\mathbb{E}\left[ \max\limits_{1\leq i\leq N}(x^i_l(0)-\bar{x}_l(0))^2\right] \right)\\
			\leq&\frac{2N}{2\lambda-\sigma^2}\sum_{l=1}^d\left(\mathbb{E}\left[ \max\limits_{1\leq i\leq N}(x^i_l(0)-\bar{x}_l(0))^2\right] \right)<\infty,\\
		\end{aligned}
	\end{equation}
	where the second inequality follows from \eqref{i}. Then, by Proposition \ref{ito}-(ii) and \eqref{mar}, we get
	\begin{equation*}
		\mathbb{E}\left[\left( \int_0^t(x^i_l(s)-\bar{x}^*_l(s))\mathrm{d}W_l(s)\right)^2\right] =\mathbb{E}\left[ \int_0^t(x^i_l(s)-\bar{x}^*_l(s))^2\mathrm{d}s\right] <\infty.
	\end{equation*}
	Therefore, $\tilde{\mathcal{I}}_l^i(t)$ is bounded in $L^1$. And by \cite[Section 3.3]{2013Brownian}, it is clear that $\tilde{\mathcal{I}}_l^i(t)$ is a martingale. Hence by Theorem \ref{marcon}, we know $\lim\limits_{t\rightarrow\infty}\tilde{\mathcal{I}}_l^i(t)$ exists a.s..
	Thus, recalling \eqref{I12}$, \lim\limits_{t\rightarrow\infty}x^i_l(t)$ exists a.s. for all $i=1,\ldots,N$ and $l=1,\ldots,d.$
	
	Now we have shown that for each $i=1,\ldots,N$, there exists a random variable $x_\infty^i$ such that
	$$\lim\limits_{t\rightarrow\infty}x^i(t)=x_\infty^i,~~{\rm a.s.}.$$
	Following from Theorem \ref{consensus1}, we have for any $1\leq i,j\leq N$,
	$$\lim\limits_{t\rightarrow\infty}\|x^i(t)-x^j(t)\|=0,~~{\rm a.s.}.$$
	Therefore, there exists a random vector $x_\infty$ such that
	\begin{equation*}
		\lim\limits_{t\rightarrow\infty}x^i(t)=x_\infty,~~{\rm a.s.},~~\forall1\leq i\leq N.
	\end{equation*}
\qed

\subsection{Error Estimation}\label{sec:4.2}

In this subsection, our concern is the error estimation of the objective function value at the state of the SCBO algorithm towards the global minimum of problem \eqref{problem}.
First of all, we give some preparations.

In what follows, for notational simplicity, we denote $\mu_t:=\mu(t)$, and for $l,k=1,\ldots,d$, 
$$\partial_l\tilde{f}(x,\mu):=\frac{\partial \tilde{f}(x,\mu)}{\partial x_l},~\partial_{kl}^2\tilde{f}(x,\mu):=\frac{\partial^2\tilde{f}(x,\mu)}{\partial x_k\partial x_l},
~\frac{\partial_\mu \tilde{f}(x,\mu_t)}{\partial t}:=\frac{\partial \tilde{f}(x,\mu_t)}{\partial \mu_t}\cdot\dot{\mu}_t,$$
where $\dot{\mu}_t=\dfrac{\partial \mu_t}{\partial t}$, and set
\begin{equation}\label{fmin}
	f_{min}:=\inf_{x\in\mathbb{R}^d}f(x).
\end{equation}
Thus, by \eqref{eq-mu3},
\begin{equation}\label{eq-L1}
	\tilde{f}(x^i(t),\mu_t)\geq f(x)-\kappa\mu^{1-q}\geq f_{min}-\kappa\bar{\mu}^{1-q}=:\tilde{f}_{min},~~\forall x\in\mathcal{X},~\mu\in(0,\bar{\mu}],
\end{equation}
where $\mathcal{X}\subseteq\mathbb{R}^d$ is a given bounded set.

First, note that
\begin{equation}\label{eq-par}
	\begin{aligned}
		\partial_k\left(\mathrm{e}^{-\beta\tilde{f}(x^i,\mu)}\right)=&-\beta \mathrm{e}^{-\beta\tilde{f}(x^i,\mu)}\partial_k\tilde{f}(x^i,\mu),\\
		\partial_{kl}^2\left(\mathrm{e}^{-\beta\tilde{f}(x^i,\mu)}\right)=&\beta \mathrm{e}^{-\beta\tilde{f}(x^i,\mu)}
		\left[\beta\partial_l\tilde{f}(x^i,\mu)\cdot\partial_k\tilde{f}(x^i,\mu)-\partial^2_{kl}\tilde{f}(x^i,\mu)\right],\\
		\frac{\partial_\mu}{\partial t}\left(\mathrm{e}^{-\beta\tilde{f}(x^i,\mu_t)}\right)=&-\beta \mathrm{e}^{-\beta\tilde{f}(x^i,\mu_t)}\nabla_\mu\tilde{f}(x^i,\mu_t)\dot{\mu}_t.
	\end{aligned}
\end{equation}
Then, applying Ito's formula to $\frac{1}{N}\sum_{i=1}^N\mathrm{e}^{-\beta\tilde{f}(x^i(t),\mu_t)}$, by Lemma \ref{lem4.3} and \eqref{eq-par}, one has
\begin{equation}\label{eq-ito}
	\resizebox{0.9\hsize}{!}{$\begin{aligned}
			&\mathrm{d}\left(\frac{1}{N}\sum_{i=1}^N\mathrm{e}^{-\beta\tilde{f}(x^i(t),\mu_t)}\right)\\
			=&\frac{1}{N}\sum_{i=1}^N\Bigg[\frac{\partial_\mu}{\partial t}\left(\mathrm{e}^{-\beta\tilde{f}(x^i(t),\mu_t)}\right)\mathrm{d}t+\sum\limits_{k=1}^d\partial_k\left(\mathrm{e}^{-\beta\tilde{f}(x^i(t),\mu_t)}\right)\mathrm{d}x^i_k(t)\\
			&+\frac{1}{2}\sum\limits_{k,l=1}^d\partial_{kl}^2\left(\mathrm{e}^{-\beta\tilde{f}(x^i(t),\mu_t)}\right)\mathrm{d}x^i_k(t)\cdot \mathrm{d}x^i_l(t)\Bigg]\\
			=:&\frac{1}{N}\mathcal{Q}_1\mathrm{d}t+\frac{1}{N}\mathcal{Q}_2\mathrm{d}t+\frac{1}{N}\mathcal{Q}_3\mathrm{d}t\\
			&-\frac{1}{N}\sum_{i=1}^N\beta\sigma \mathrm{e}^{-\beta\tilde{f}(x^i(t),\mu_t)}\nabla_x\tilde{f}(x^i(t),\mu_t)^\mathrm{T}
			\left[\sum\limits_{k=1}^d(x^i_k(t)-\bar{x}^*_k(t))\mathrm{d}W_k(t)e_k\right],
		\end{aligned}$}
\end{equation}
where
\begin{equation*}
	\small
	\begin{aligned}
		\mathcal{Q}_1=&-\sum_{i=1}^N\left( \beta \mathrm{e}^{-\beta\tilde{f}(x^i(t),\mu_t)}\nabla_\mu\tilde{f}(x^i(t),\mu_t)\dot{\mu}_t\right) ,\\
		\mathcal{Q}_2=&\sum_{i=1}^N\left( \beta\lambda \mathrm{e}^{-\beta\tilde{f}(x^i(t),\mu_t)}\nabla_x\tilde{f}(x^i(t),\mu_t)^\mathrm{T}(x^i(t)-\bar{x}^*(t))\right) ,\\
		\mathcal{Q}_3=& \sum_{i=1}^N\left[\frac{\beta\sigma^2\mathrm{e}^{-\beta\tilde{f}(x^i(t),\mu_t)}}{2} \sum\limits_{k=1}^d\left(-\partial^2_{kk}\tilde{f}(x^i(t),\mu_t)
		+\beta(\partial_k\tilde{f}(x^i(t),\mu_t))^2\right)
		(x^i_k(t)-\bar{x}^*_k(t))^2\right].\\
	\end{aligned}
\end{equation*}
Taking expectation on the both sides of (\ref{eq-ito}) and using Proposition \ref{ito}-(i), we get
\begin{equation}\label{eq-exp}
	\mathrm{d}\left(\frac{1}{N}\sum_{i=1}^N\mathbb{E}\left[ \mathrm{e}^{-\beta\tilde{f}(x^i(t),\mu_t)}\right] \right)=\frac{1}{N}\left( \mathbb{E}\left[\mathcal{Q}_1\right] \mathrm{d}t+\mathbb{E}\left[\mathcal{Q}_2\right] \mathrm{d}t+\mathbb{E}\left[\mathcal{Q}_3\right] \mathrm{d}t\right) .
\end{equation}

Next, we separately give the estimations on the terms of $\mathbb{E}\left[ \mathcal{Q}_i\right] ,~i=1,2,3$.

$\bullet$ (Estimate on $\mathbb{E}\left[ \mathcal{Q}_1\right]$): According to the proof of Theorem \ref{consensus} and \eqref{alg-2}, there exists a bounded set $\mathcal{X}$ such that for all $i=1,\ldots,N$,
\begin{equation}\label{X}
	x^i(t)\in\mathcal{X}\quad\mbox{and}\quad \bar{x}^*(t)\in\mathcal{X},~~{\rm a.s.},~~\forall t\geq 0.
\end{equation}
By \eqref{eq-mu1}, \eqref{eq-L1}, \eqref{X} and the monotone decreasing of $\mu(\cdot)$ on $(0,+\infty)$, we have
\begin{equation}\label{eq-I1}
	\mathbb{E}\left[\mathcal{Q}_1\right] 
	\geq\sum_{i=1}^N\beta\kappa \mathrm{e}^{-\beta\tilde{f}_{min}}\mu_t^{-q}\dot{\mu}_t
	=N\beta\kappa \mathrm{e}^{-\beta\tilde{f}_{min}}\mu_t^{-q}\dot{\mu}_t.
\end{equation}

$\bullet$ (Estimate on $\mathbb{E}\left[ \mathcal{Q}_2\right]$): First, by the definition of $\bar{x}^*(t)$ in \eqref{alg-2}, we have
$$\left(\sum\limits_{i=1}^N\mathrm{e}^{-\beta\tilde{f}(x^i(t),\mu_t)}\right)\bar{x}^*(t)=\sum\limits_{i=1}^N\mathrm{e}^{-\beta\tilde{f}(x^i(t),\mu_t)}x^i(t).$$
This yields
\begin{equation}\label{eq-0}
	\sum\limits_{i=1}^N\mathrm{e}^{-\beta\tilde{f}(x^i(t),\mu_t)}\nabla_x\tilde{f}(\bar{x}^*(t),\mu_t)^\mathrm{T}(\bar{x}^*(t)-x^i(t))=0.
\end{equation}
Then, using \eqref{eq-mu1}, \eqref{eq-L1}, \eqref{X} and \eqref{eq-0}, it shows that
\begin{equation}\label{eq-I2}
	\resizebox{0.9\hsize}{!}{$\begin{aligned}
			&\mathbb{E}\left[\mathcal{Q}_2\right] \\
			=&\beta\lambda\sum\limits_{i=1}^N\mathbb{E}\left[\mathrm{e}^{-\beta\tilde{f}(x^i(t),\mu_t)}\left(\nabla_x\tilde{f}(x^i(t),\mu_t)-\nabla_x\tilde{f}(\bar{x}^*(t),\mu_t)\right)^\mathrm{T}(x^i(t)-\bar{x}^*(t))\right]\\
			\geq&-\beta\lambda\sum\limits_{i=1}^N\mathbb{E}\left[\mathrm{e}^{-\beta\tilde{f}(x^i(t),\mu_t)}\sup\limits_{(x,\mu)\in\mathcal{X}\times(0,\bar{\mu}]}\left\|\nabla_{xx}^2\tilde{f}(x,\mu)\right\|_2\| x^i(t)-\bar{x}^*(t)\|^2\right]\\
			\geq&-\lambda \eta\beta \mathrm{e}^{-\beta\tilde{f}_{min}}\mu_t^{-q-1}\sum\limits_{i=1}^N\mathbb{E}\left[ \|x^i(t)-\bar{x}^*(t)\|^2\right] .\\
		\end{aligned}$}
\end{equation}

$\bullet$ (Estimate on $\mathbb{E}\left[ \mathcal{Q}_3\right]$): By \eqref{eq-mu1}, one has
$$\max\limits_{1\leq k\leq d}\sup\limits_{(x,\mu)\in\mathcal{X}\times(0,\bar{\mu}]}\left|\partial^2_{kk}\tilde{f}(x,\mu)\right|\leq\sup\limits_{(x,\mu)\in\mathcal{X}\times(0,\bar{\mu}]}\left\|\nabla_{xx}^2\tilde{f}(x,\mu)\right\|_2 \leq \eta\mu^{-q-1}.$$
Then, it follows from \eqref{eq-L1} and \eqref{X} that
\begin{equation}\label{eq-I3}
	\begin{aligned}
		\mathbb{E}\left[\mathcal{Q}_3\right] &\geq-\sigma^2\beta\sum\limits_{i=1}^N\mathbb{E}\left[\mathrm{e}^{-\beta\tilde{f}(x^i(t),\mu_t)}\sum\limits_{k=1}^d\partial^2_{kk}\tilde{f}(x^i(t),\mu_t)(x^i_k(t)-\bar{x}^*_k(t))^2\right]\\
		&\geq-\sigma^2\eta\beta \mathrm{e}^{-\beta\tilde{f}_{min}}\mu_t^{-q-1}\sum_{i=1}^N\mathbb{E}\left[ \|x^i(t)-\bar{x}^*(t)\|^2\right] .\\
	\end{aligned}
\end{equation}
Before giving the error analysis, we first list some necessary assumptions.
\begin{assumption}\label{ass1}
	System parameters $\lambda$ and $\sigma$ in the SCBO algorithm satisfy $2\lambda>\sigma^2$.
\end{assumption}
\begin{assumption}\label{ass2}
	Initial data $\{x^i(0)\}_{i=1}^N$ is i.i.d and $x^i(0)\sim x^{\mathrm{in}},$ where $x^{\mathrm{in}}$ is a reference random variable.
\end{assumption}
\begin{assumption}\label{ass3}
	There exist constants $T>0$ and $\gamma>0$ such that function $\mu:\mathbb{R}_+\rightarrow(0,\bar{\mu}]$ satisfies
	\begin{equation}\label{condition2}
		\int_0^t\mathrm{e}^{-(2\lambda-\sigma^2)s}(\mu(s))^{-q-1}\mathrm{d}s\leq \gamma, ~~\forall t\geq T.
	\end{equation}
\end{assumption}
\begin{remark}
	Under Assumption \ref{ass1}, Assumption \ref{ass3} is much trival. For example, if we set $\mu_t=\mathrm{e}^{-\alpha t}$ with $0<(q+1)\alpha<2\lambda-\sigma^2$, then Assumption \ref{ass3} holds.
\end{remark}

\begin{theorem}\label{th-con}
	Assume that Assumptions \ref{ass1}-\ref{ass3} hold. If the initial data $\left\lbrace \mu_0,x^i(0)\right\rbrace_{i=1}^N$ and system parameters $\{\beta, \lambda, \sigma\}$ in the SCBO algorithm satisfy
	\begin{equation}\label{condition}
		\resizebox{0.9\hsize}{!}{$\begin{aligned}
				&\left(\mathrm{e}^{-\mu_0^{1-q}}\beta\kappa-\varepsilon\right)\mathbb{E}\left[\mathrm{e}^{-\beta f(x^{\mathrm{in}})}\right]\\
				\geq &\frac{\mu_0^{1-q}}{1-q}\beta\kappa \mathrm{e}^{-\beta\tilde{f}_{min}}+\gamma\left(2\lambda+\sigma^2\right)\eta\beta \mathrm{e}^{-\beta\tilde{f}_{min}}
				\sum\limits_{l=1}^d\mathbb{E}\left[\max\limits_{1\leq i\leq N}(x^i_l(0)-\bar{x}_l(0))^2\right],
			\end{aligned}$}
	\end{equation}
	for some $0<\varepsilon<\mathrm{e}^{-\mu_0^{1-q}\beta\kappa}$. Then, one has
	$${\rm ess}\inf f(x_\infty)\leq f_{min}+E(\beta),$$
	where $x_\infty$ is the consensus state defined as in Theorem \ref{consensus}, and $\lim_{\beta\rightarrow\infty}E(\beta)=0$.
\end{theorem}
{\noindent\it Proof}
	By (\ref{eq-exp}), we have
	\begin{equation}\label{eq-p0}
		\resizebox{0.9\hsize}{!}{$\begin{aligned}
				&\frac{\mathrm{d}}{\mathrm{d}t}\left(\frac{1}{N}\sum_{i=1}^N\mathbb{E}\left[\mathrm{e}^{-\beta\tilde{f}(x^i(t),\mu_t)}\right]\right) \\
				\geq&\beta\kappa \mathrm{e}^{-\beta\tilde{f}_{min}}\mu_t^{-q}\dot{\mu}_t-\left(\lambda+\frac{1}{2}\sigma^2\right)\eta\beta \mathrm{e}^{-\beta\tilde{f}_{min}}\mu_t^{-q-1}\frac{1}{N}\sum\limits_{i=1}^N\mathbb{E}\left[ \|x^i(t)-\bar{x}^*(t)\|^2\right] \\
				\geq&\beta\kappa \mathrm{e}^{-\beta\tilde{f}_{min}}\mu_t^{-q}\dot{\mu}_t
				-\frac{(2\lambda+\sigma^2)\eta\beta \mathrm{e}^{-\beta\tilde{f}_{min}}\mathrm{e}^{-(2\lambda-\sigma^2)t}}{\mu_t^{q+1}}
				\sum\limits_{l=1}^d\mathbb{E}\left[\max\limits_{1\leq i\leq N}(x^i_l(0)-\bar{x}_l(0))^2\right],\\
			\end{aligned}$}
	\end{equation}
	where the first inequality follows from \eqref{eq-I1}, \eqref{eq-I2} and \eqref{eq-I3}, and the last inequality follows from Lemma \ref{3}-(iii).
	Now, integrating (\ref{eq-p0}) in $t$ gives
	\begin{equation}\label{eq-p1}
		\resizebox{0.9\hsize}{!}{$\begin{aligned}
				&\frac{1}{N}\sum_{i=1}^N\mathbb{E}\left[ \mathrm{e}^{-\beta\tilde{f}(x^i(t),\mu_t)}\right] \\
				\geq&\frac{1}{N}\sum_{i=1}^N\mathbb{E}\left[ \mathrm{e}^{-\beta\tilde{f}(x^i(0),\mu_0)}\right] +\int_0^t\left(\beta\kappa \mathrm{e}^{-\beta\tilde{f}_{min}}\mu_s^{-q}\dot{\mu}_s\right)\mathrm{d}s\\
				&-(2\lambda+\sigma^2)\eta\beta \mathrm{e}^{-\beta\tilde{f}_{min}}\sum\limits_{l=1}^d\mathbb{E}\left[\max\limits_{1\leq i\leq N}(x^i_l(0)-\bar{x}_l(0))^2\right]\left(\int_0^t\mathrm{e}^{-(2\lambda-\sigma^2)s}\mu_s^{-q-1}\mathrm{d}s\right).\\
			\end{aligned}$}
	\end{equation}
	By $0\leq q <1$, we have
	\begin{equation}\label{eq-p2}
		\int_0^t\left(\beta\kappa \mathrm{e}^{-\beta\tilde{f}_{min}}\mu_s^{-q}\dot{\mu}_s\right)\mathrm{d}s=\frac{1}{1-q}\beta\kappa \mathrm{e}^{-\beta\tilde{f}_{min}}\left(\mu_t^{1-q}-\mu_0^{1-q}\right).
	\end{equation}
	By \eqref{eq-mu3} and \eqref{X}, we get
	\begin{equation}\label{eq-p4}
		\mathbb{E}\left[ \mathrm{e}^{-\beta\tilde{f}(x^i(0),\mu_0)}\right] \geq\mathbb{E}\left[\mathrm{e}^{-\beta(f(x^i(0))+\kappa\mu_0^{1-q})}\right] 
		=\mathrm{e}^{-\beta\kappa\mu_0^{1-q}}\mathbb{E}\left[\mathrm{e}^{-\beta f(x^i(0))}\right] .
	\end{equation}
	Then, substituting \eqref{eq-p2} and \eqref{eq-p4} into \eqref{eq-p1}, one has
	\begin{equation}\label{eq-p5}
		\resizebox{0.9\hsize}{!}{$\begin{aligned}
				&\frac{1}{N}\sum_{i=1}^N\mathbb{E}\left[\mathrm{e}^{-\beta\tilde{f}(x^i(t),\mu_t)}\right] \\
				\geq&\frac{1}{N}\sum_{i=1}^N\mathrm{e}^{-\beta\kappa\mu_0^{1-q}}\mathbb{E}\left[\mathrm{e}^{-\beta f(x^i(0))}\right]
				+\frac{1}{1-q}\beta\kappa \mathrm{e}^{-\beta\tilde{f}_{min}}\left(\mu_t^{1-q}-\mu_0^{1-q}\right)\\
				&-(2\lambda+\sigma^2)\eta\beta \mathrm{e}^{-\beta\tilde{f}_{min}}\sum\limits_{l=1}^d\mathbb{E}\left[\max\limits_{1\leq i\leq N}(x^i_l(0)-\bar{x}_l(0))^2\right]\left(\int_0^t\mathrm{e}^{-(2\lambda-\sigma^2)s}\mu_s^{-q-1}\mathrm{d}s\right).\\
			\end{aligned}$}
	\end{equation}
	Following from Definition \ref{def-sf}-(ii), we have
	$$\lim_{t\rightarrow\infty}\tilde{f}(x^i(t),\mu_t)=f(x_\infty),~~{\rm a.s.},~~\forall i=1,2,\ldots,N.$$
	This yields that
	\begin{equation}\label{limit}
		\lim_{t\rightarrow\infty}\frac{1}{N}\sum_{i=1}^N\mathbb{E}\left[\mathrm{e}^{-\beta\tilde{f}(x^i(t),\mu_t)}\right]=\mathbb{E}\left[\mathrm{e}^{-\beta f(x_\infty)}\right].
	\end{equation}
	Letting $t\rightarrow\infty$ in \eqref{eq-p5}, and using  \eqref{limit}, it has
	\begin{equation}\label{eq-p6}
		\begin{aligned}
			&\mathbb{E}\left[\mathrm{e}^{-\beta f(x_\infty)}\right]\\
			\geq&\frac{1}{N}\sum_{i=1}^N\mathrm{e}^{-\beta\kappa\mu_0^{1-q}}\mathbb{E}\left[\mathrm{e}^{-\beta f(x^{\mathrm{in}})}\right]-\frac{\mu_0^{1-q}}{1-q}\beta\kappa \mathrm{e}^{-\beta\tilde{f}_{min}}\\
			&-\gamma\left(2\lambda+\sigma^2\right)\eta\beta \mathrm{e}^{-\beta\tilde{f}_{min}}
			\sum\limits_{l=1}^d\mathbb{E}\Big[\max\limits_{1\leq i\leq N}(x^i_l(0)-\bar{x}_l(0))^2\Big]\\
			\geq&\varepsilon\mathbb{E}\mathrm{e}^{-\beta f(x^{\mathrm{in}})},
		\end{aligned}
	\end{equation}
	where the first inequality uses Assumptions \ref{ass2} and \ref{ass3}, and the last inequality
	follows from \eqref{condition}.
	Therefore, from \eqref{eq-p6} we have
	\begin{equation}\label{log}
		\mathrm{e}^{-\beta{\rm ess}\inf f(x_\infty)}=\mathbb{E}\left[\mathrm{e}^{-\beta{\rm ess}\inf f(x_\infty)}\right]\geq\mathbb{E}\left[\mathrm{e}^{-\beta f(x_\infty)}\right]\geq\varepsilon\mathbb{E}\left[\mathrm{e}^{-\beta f(x^{\mathrm{in}})}\right].
	\end{equation}
	Take logarithm and divide by $-\beta$ to the both sides of \eqref{log} to find
	\begin{equation*}
		\begin{aligned}
			{\rm ess}\inf f(x_\infty)\leq& -\frac{1}{\beta}\log\left(\varepsilon\mathbb{E}\left[\mathrm{e}^{-\beta f(x^{\mathrm{in}})}\right]\right)
			=-\frac{1}{\beta}\log\mathbb{E}\left[\mathrm{e}^{-\beta f(x^{\mathrm{in}})}\right]-\frac{1}{\beta}\log\varepsilon.
		\end{aligned}
	\end{equation*}
	So if we define
	\begin{equation}\label{Ebetaa}
		E(\beta):=-\frac{1}{\beta}\log\mathbb{E}\left[\mathrm{e}^{-\beta f(x^{\mathrm{in}})}\right] -f_{min}-\frac{1}{\beta}\log\varepsilon,
	\end{equation}
	we have
	$${\rm ess}\inf f(x_\infty)\leq f_{min}+E(\beta),$$
	where $\lim_{\beta\rightarrow\infty}E(\beta)=0$ by Proposition \ref{Laplace}.
\qed
\begin{remark}
	By Remark \ref{q}, if $f$ is locally Lipschitz continuous and we can construct a smoothing function of it with $q=0$, then \eqref{condition} is
	\begin{equation*}
		\begin{aligned}
			&\left(\mathrm{e}^{-\mu_0\beta\kappa}-\varepsilon\right)\mathbb{E}\left[\mathrm{e}^{-\beta f(x^{\mathrm{in}})}\right]\\
			\geq &\mu_0\beta\kappa \mathrm{e}^{-\beta\tilde{f}_{min}}+\gamma\left(2\lambda+\sigma^2\right)\eta\beta \mathrm{e}^{-\beta\tilde{f}_{min}}
			\sum\limits_{l=1}^d\mathbb{E}\left[\max\limits_{1\leq i\leq N}(x^i_l(0)-\bar{x}_l(0))^2\right].
		\end{aligned}
	\end{equation*}
	If $f$ is defined by $f:=|g(x)|^p$ as in Remark \ref{q}, then we can give a smoothing function of it with $q=1-p$, and \eqref{condition} can be rewritten into
	\begin{equation*}
		\begin{aligned}
			&\left(\mathrm{e}^{-\mu_0^{p}\beta\kappa}-\varepsilon\right)\mathbb{E}\left[\mathrm{e}^{-\beta f(x^{\mathrm{in}})}\right]\\
			\geq &\frac{\mu_0^{p}}{p}\beta\kappa \mathrm{e}^{-\beta\tilde{f}_{min}}+\gamma\left(2\lambda+\sigma^2\right)\eta\beta \mathrm{e}^{-\beta\tilde{f}_{min}}
			\sum\limits_{l=1}^d\mathbb{E}\left[\max\limits_{1\leq i\leq N}(x^i_l(0)-\bar{x}_l(0))^2\right].
		\end{aligned}
	\end{equation*}
\end{remark}

\begin{remark}
	By multiplying $\mathrm{e}^{\beta\tilde{f}_{min}}$ on the both sides of \eqref{condition}, we obtain
	\begin{equation}\label{beta}
		\begin{aligned}
			&\left(\mathrm{e}^{-\mu_0^{1-q}\beta\kappa}-\varepsilon\right)\mathbb{E}\left[\mathrm{e}^{\beta\left(\tilde{f}_{min}- f(x^{\mathrm{in}})\right)}\right]\\
			\geq& \frac{\mu_0^{1-q}}{1-q}\beta\kappa +\gamma\left(2\lambda+\sigma^2\right)\eta\beta\sum\limits_{l=1}^d\mathbb{E}\left[\max\limits_{1\leq i\leq N}(x^i_l(0)-\bar{x}_l(0))^2\right].
		\end{aligned}
	\end{equation}
	If let $\beta\rightarrow0$ in \eqref{beta}, we can see that the left side of \eqref{beta} converges to $1-\varepsilon$, and the right side converges to $0$.
	Therefore, $\beta>0$ satisfying \eqref{condition} is guaranteed to exist.
	But if let $\beta\rightarrow\infty$, by $\tilde{f}_{min}- f(x^{\mathrm{in}})<0$, we can see that the left side of \eqref{beta} is not larger than $0$
	and the right side of it tends to $\infty$. Hence, the estimate in Theorem \ref{th-con} is only an error estimate for the SCBO algorithm with a proper $\beta$, but can not guarantee a sufficient small error estimate on the objective function value to $f_{min}$. A similar estimation was given in \cite{HaConvergence2020} to \eqref{problem} with $f\in\mathcal{C}_b^2\left(\mathbb{R}^d \right) $.
	
	In what follows, we will give a sufficient condition to guarantee the small enough of the error $E(\beta)$.
	It illustrates that when the parameters are proper selected, the objective function value at the global consensus state can be close enough to the global minimum. It is noteworthy that no previous literature, to our knowledge, has investigated this thing.
\end{remark}
\begin{corollary}\label{co-con}
	Assume that Assumptions \ref{ass1}-\ref{ass3} hold. Then, for any given $\delta>0$, there exist $\varepsilon$ and system parameters $\{\beta, \mu_0, \lambda, \sigma\}$ in the SCBO algorithm satisfying \eqref{condition}. Moreover, with these parameters, one has
	$${\rm ess}\inf f(x_\infty)\leq f_{min}+E(\beta),$$
	where $x_\infty$ is the consensus state of the SCBO algorithm defined as in Theorem \ref{consensus}, and $E(\beta)$ defined in \eqref{beta} satisfies $E(\beta)\leq \delta$.
\end{corollary}
{\noindent\it Proof}
	For any given $\delta>0$, let $\varepsilon= \frac{\mathrm{e}^{-2\mu_0^{1-q}\beta\kappa-\beta\delta}}{\frac{1}{1-q}\mu_0^{1-q}\beta\kappa+\mathrm{e}^{-\mu_0^{1-q}\beta\kappa-\beta\delta}}$. In what follows, we will show that there exist corresponding parameters that make them satisfy \eqref{condition}.
	
	Firstly, for the given $\delta>0$ and $\beta>0$, there exists a $\mu_0>0$ satisfying
	\begin{equation}\label{re3}
		\delta>-\frac{1}{\beta}\ln\left( 1-\frac{\mu_0^{1-q}\beta\kappa}{1-q}\mathrm{e}^{\mu_0^{1-q}\beta\kappa}\right),
	\end{equation}
	in which if we let $\mu_0\rightarrow 0$ in \eqref{re3}, then we can see that the right side of it convergences to $0$ and the left side equals $\delta>0$. By some transformation on \eqref{re3}, one has
	\begin{equation}\label{re2}
		\mathrm{e}^{-\mu_0^{1-q}\beta\kappa}(1-\mathrm{e}^{-\beta\delta})-\frac{1}{1-q}\mu_0^{1-q}\beta\kappa>0.
	\end{equation}
	
	Next, considering
	\begin{equation*}
		\mathrm{e}^{-\mu_0^{1-q}\beta\kappa}-\varepsilon-\frac{1}{1-q}\mu_0^{1-q}\beta\kappa=\frac{\frac{1}{1-q}\mu_0^{1-q}\beta\kappa\left(\mathrm{e}^{-\mu_0^{1-q}\beta\kappa}(1-\mathrm{e}^{-\beta\delta})-\frac{1}{1-q}\mu_0^{1-q}\beta\kappa \right) }{\frac{1}{1-q}\mu_0^{1-q}\beta\kappa+\mathrm{e}^{-\mu_0^{1-q}\beta\kappa-\beta\delta}}
	\end{equation*}
	and $0\leq q<1$, we know that for the given $\delta>0$ and $\beta>0$, there exists a sufficiently small $\mu_0>0$ such that
	\begin{equation}\label{re1}
		\left(\mathrm{e}^{-\mu_0^{1-q}\beta\kappa}-\varepsilon\right)-\frac{1}{1-q}\mu_0^{1-q}\beta\kappa>0.
	\end{equation} 
	Consequently, following from \eqref{re1}, we can find a sufficiently small $\beta>0$ and a suitable $x^{\mathrm{in}}$ satisfying
	\begin{equation}\label{re4}
		\left(\mathrm{e}^{-\mu_0^{1-q}\beta\kappa}-\varepsilon\right)\mathbb{E}\left[\mathrm{e}^{\beta\left(\tilde{f}_{min}- f(x^{\mathrm{in}})\right)}\right]-\frac{1}{1-q}\mu_0^{1-q}\beta\kappa>0.
	\end{equation}
	
	At last, from Assumption \ref{ass3}, we know $\gamma$ depends on $\lambda$ and $\sigma$.
	In view of \eqref{re4}, there are appropriate $\lambda$ and $\sigma$ to ensure
	$$(2\lambda+\sigma^2)\gamma\leq\frac{\left(\mathrm{e}^{-\mu_0^{1-q}\beta\kappa}-\varepsilon\right)\mathbb{E}\left[\mathrm{e}^{\beta\left(\tilde{f}_{min}- f(x^{\mathrm{in}})\right)}\right]-\frac{1}{1-q}\mu_0^{1-q}\beta\kappa}
	{\eta\beta \sum\limits_{l=1}^d\mathbb{E}\left[\max\limits_{1\leq i\leq N}(x^i_l(0)-\bar{x}_l(0))^2\right]},$$
	which suggests \eqref{beta} as well as \eqref{condition} holds.
	
	Putting all the above results together, we confirm that for any $\delta>0$, there exist suitable parameters $\varepsilon$, $\beta$, $\mu_0$, $\lambda$ and $\sigma$ such that \eqref{condition} holds.
	
	In addition, from Theorem \ref{th-con}, we have
	$${\rm ess}\inf f(x_\infty)\leq f_{min}+E(\beta),$$
	where $E(\beta)$ is defined in \eqref{Ebetaa}.
	In order to express conveniently, we define
	$$\mathcal{L}:=\frac{\left(\frac{1}{1-q}\mu_0^{1-q}\kappa+\gamma\left(2\lambda+\sigma^2\right)\eta\sum\limits_{l=1}^d\mathbb{E}\left[\max\limits_{1\leq i\leq N}(x^i_l(0)-\bar{x}_l(0))^2\right]\right)\beta\mathrm{e}^{-\beta\tilde{f}_{min}}}{\mathrm{e}^{-\mu_0^{1-q}\beta\kappa}-\varepsilon},$$
	where $\{\varepsilon, \beta, \mu_0, \lambda, \sigma \}$  is the aforementioned parameters that make \eqref{condition} valid.
	
	On the one hand, \eqref{condition} implies that
	\begin{equation}\label{Ebeta}
		E(\beta)\leq-\frac{1}{\beta}\log\mathcal{L}-f_{min}-\frac{1}{\beta}\log\varepsilon.
	\end{equation}
	
	On the other hand, by $\varepsilon= \frac{\mathrm{e}^{-2\mu_0^{1-q}\beta\kappa-\beta\delta}}{\frac{1}{1-q}\mu_0^{1-q}\beta\kappa+\mathrm{e}^{-\mu_0^{1-q}\beta\kappa-\beta\delta}}$, we have
	$$\frac{\frac{1}{1-q}\mu_0^{1-q}\beta\kappa\varepsilon}{\mathrm{e}^{-\mu_0^{1-q}\beta\kappa}-\varepsilon}= \mathrm{e}^{-\beta(\delta+\kappa\mu_0^{1-q})}.$$
	Then, we have
	\begin{equation*}
		\frac{\left(\frac{1}{1-q}\mu_0^{1-q}\kappa+\gamma\left(2\lambda+\sigma^2\right)\eta\sum\limits_{l=1}^d\mathbb{E}\left[\max\limits_{1\leq i\leq N}(x^i_l(0)-\bar{x}_l(0))^2\right]\right)\beta\varepsilon}{\mathrm{e}^{-\mu_0^{1-q}\beta\kappa}-\varepsilon}\geq\mathrm{e}^{-\beta(\delta+\kappa\mu_0^{1-q})}.
	\end{equation*}
	We multiply the both sides of the above inequality by $\mathrm{e}^{-\beta\tilde{f}_{min}}$ to obtain
	\begin{equation*}
		\mathcal{L}\varepsilon\geq\mathrm{e}^{-\beta(\delta+\kappa\mu_0^{1-q}+\tilde{f}_{min})}.
	\end{equation*}
	The definition of $\tilde{f}_{min}$ in \eqref{eq-L1} along with the fact that $\mu_0\leq\bar{\mu}$ yields
	\begin{equation}\label{delta2}
		\mathcal{L}\varepsilon\geq\mathrm{e}^{-\beta(\delta+f_{min})}.
	\end{equation}
	Taking logarithm and dividing by $-\beta$ to the both sides of \eqref{delta2}, we have
	\begin{equation*}
		-\frac{1}{\beta}\log\left(\mathcal{L}\varepsilon\right) \leq\delta+f_{min},
	\end{equation*}
	i.e.
	\begin{equation}\label{delta3}
		-\frac{1}{\beta}\log\mathcal{L}-f_{min}-\frac{1}{\beta}\log\varepsilon\leq\delta.
	\end{equation}
	Therefore, combining \eqref{Ebeta} and \eqref{delta3}, we obtain $E(\beta)\leq \delta$.
\qed

\section{Numerical Simulations}
\label{sec:5}

In the section that follows, we give two numerical examples to verify our theoretical results and illustrate the good performance of the proposed SCBO algorithm.
All the numerical testings are performed on a Lenovo PC using PYTHON 3.8.10.
For numerical simulations, inspired by the proposed predictor-corrector type model in \cite{HaConvergence2021}, we use the two-step numerical scheme of the SCBO algorithm with $h=\Delta t>0$ in Algorithm \ref{discrete} to realize the SCBO algorithm in all experiments, where we set $h=0.01$, let $\{w_l\}_{l=1}^d$ be i.i.d and follow the standard normal distribution. And we update the smoothing parameters by $\mu_t=\mu_0\mathrm{e}^{-\alpha t}$ with $\mu_0>0$ and $\alpha>0$.
\begin{algorithm}
	\normalsize
	\caption{D-SCBO algorithm}
	\label{discrete}
	\begin{algorithmic}
		\STATE{\begin{flalign*}
				\begin{split}
					&\hat{x}^i(nh)=\bar{x}^*(nh)+\mathrm{e}^{-\lambda h}(x^i(nh)-\bar{x}^*(nh)).\\
					&x^i((n+1)h)=\hat{x}^i(nh)-\sum_{l=1}^d(\hat{x}^i_l(nh)-\bar{x}^*_l(nh))\sigma\sqrt{h}w_l(nh)e_l.\\
					&\bar{x}^*(nh):=\dfrac{\sum_{i=1}^Nx^i(nh)\mathrm{e}^{-\beta \tilde{f}(x^i(nh),\mu(nh))}}{\sum_{i=1}^N\mathrm{e}^{-\beta \tilde{f}(x^i(nh),\mu(nh))}}.
				\end{split}&
		\end{flalign*}}
	\end{algorithmic}
\end{algorithm}

\begin{example}\label{example1}
	In this example, we test the theoretical results of the SCBO algorithm by minimizing the following nonconvex and nonsmooth objective function:
	\begin{equation}\label{examplefunction1}
		f(x)=\frac{1}{10}\sum_{l=1}^d[|x_l|-\cos(\pi x_l)+1],
	\end{equation}
	where $x=(x_1,\ldots,x_d)^\mathrm{T}\in\mathbb{R}^d$. 
	It is clear that $f$ has many local minimizers and a unique global minimizer $x^*=0\in\mathbb{R}^d$, with global optimal value $f_{min}=0$.
	
	For the nonsmooth term $|x_l|$ in \eqref{examplefunction1}, its smoothing function can be chosen as $\varphi_1(x_l,\mu)$ in \eqref{smoothing}.
	Then, a smoothing function of \eqref{examplefunction1} with $q=0$ can be defined as follows:
	\begin{equation*}
		\tilde{f}(x,\mu)=\frac{1}{10}\sum_{l=1}^d[\varphi_1(x_l,\mu)-\cos(\pi x_l)+1].
	\end{equation*}
	
	In what follows, we conduct two numerical experiments. 
	
	\textbf{Experiment 1}: Test the existence of the common consensus state of the SCBO algorithm. Set the dimension $d=3$ and the number of particles $N=100$. The initial data is uniformly chosen from the hypercube $[-2,2]^3$, and the parameters in the SCBO algorithm are set by:
	$$\beta=10,\quad \lambda=\sigma=1,\quad \mu_0=0.5,\quad \alpha=0.1.$$
	\begin{figure}[h]
		\centering
		\begin{minipage}{0.47\linewidth}
			\vspace{-2pt}
			\centerline{\includegraphics[width=\textwidth]{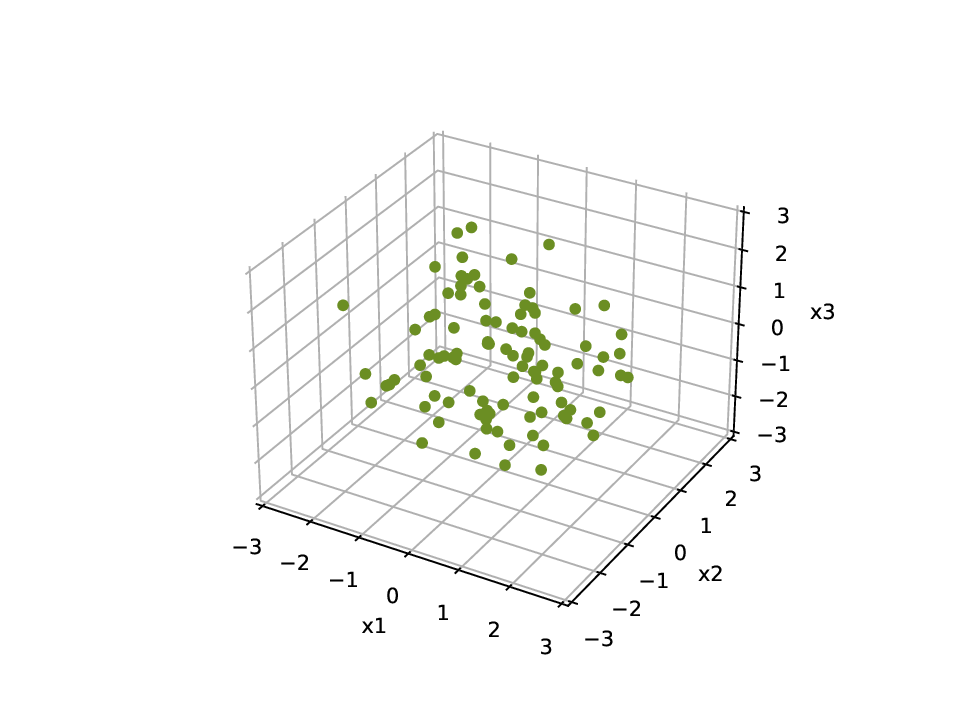}}
			\centering{$(a)$ $t=0$}
			\vspace{-2pt}
			\centerline{\includegraphics[width=\textwidth]{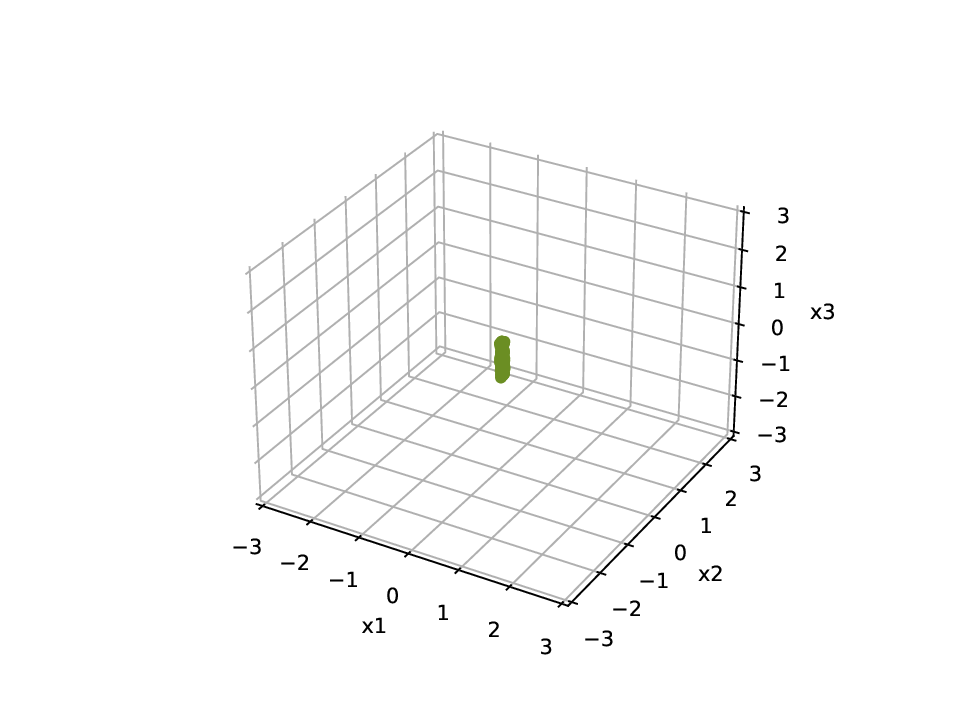}}
			\centering{$(c)$ $t=2$}
		\end{minipage}
		\hspace{-0.45in}
		\begin{minipage}{0.47\linewidth}
			\vspace{-2pt}
			\centerline{\includegraphics[width=\textwidth]{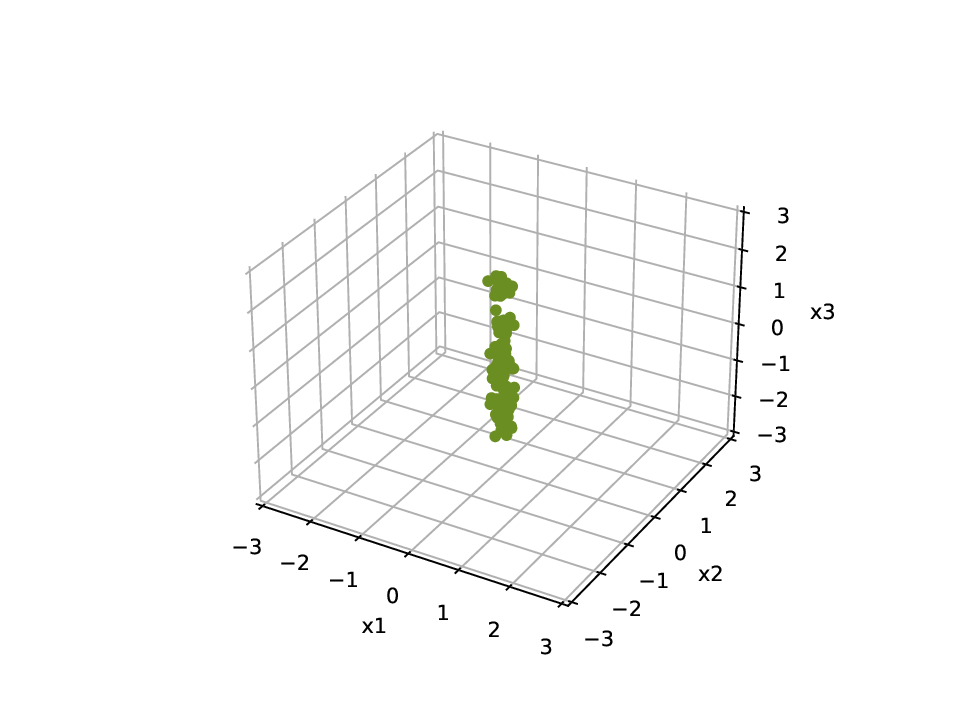}}
			\centering{$(b)$ $t=1$}
			\vspace{-2pt}
			\centerline{\includegraphics[width=\textwidth]{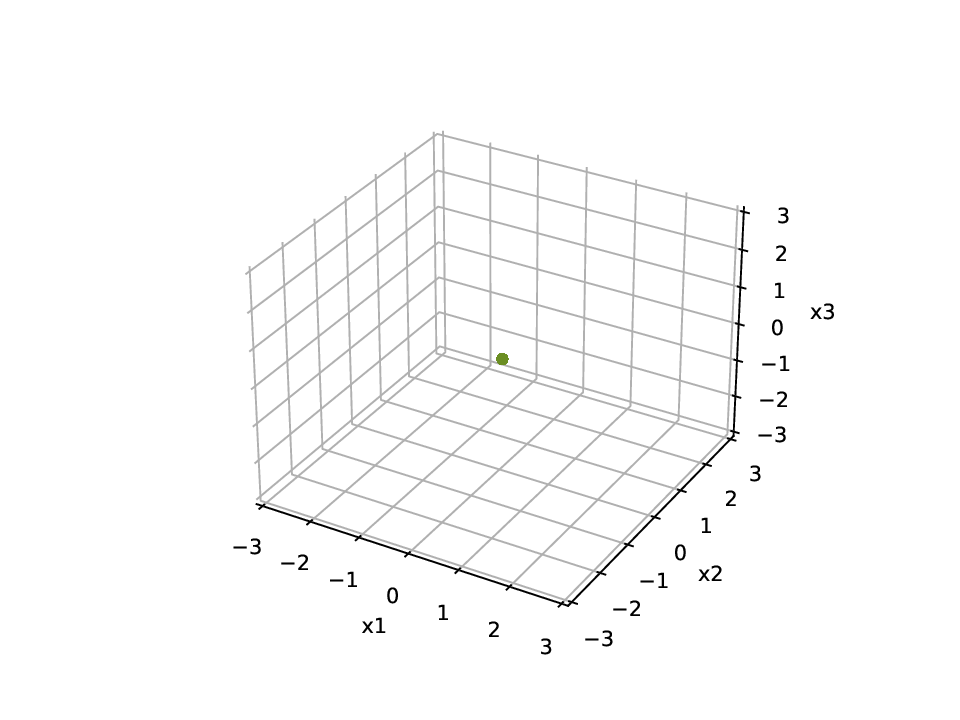}}
			\centering{$(d)$ $t=7$}
		\end{minipage}
		\caption{Evolution of the state of $100$ particles for times $t=0,1,2,7$.}
		\label{fig2}
	\end{figure}
	Fig. \ref{fig2} presents the temporal evolution of particle state for $t=0,1,2,7$ respectively.
	It can be seen that all the $100$ particles reach a common consensus state when $t=7$.
	This experiment supports the result that the SCBO algorithm exhibits a global consensus and reaches a common consensus state ultimately.
	
	\textbf{Experiment 2}: Demonstrate the error estimate between the objective function value at the common consensus state and the global minimum,
	which is used to verify the results in Corollary \ref{co-con} with $q=0$.
	Set the dimension $d=1$, the number of particles $N=150$, and choose the initial data uniformly distributed in $[-0.2,0.2]$.
	For each particle, we randomly generate $150$ initial values to get the expectation values in \eqref{condition}.
	The upper bound of error $E(\beta)$ is set by $\delta=0.01$. In order to satisfy the conditions in Corollary \ref{co-con}, we set the parameters as follows:
	$$\beta=0.2,\quad \lambda=0.1,\quad \sigma=0.3,\quad \mu_0=0.0005,\quad \alpha=0.1.$$
	The results obtained by the SCBO algorithm are given in Fig. \ref{fig3}, which presents the state trajectories of $150$ particles in the experiment.
	Note that the global consensus emerges and $x_\infty=3.19\mathrm{e}$-$04$, which is a good approximation of the minimizer $x^*=0$. Hence by calculation,
	we find that $f(x_\infty)=3.20\mathrm{e}$-$05$ and $E(\beta)=1.09\mathrm{e}$-$04$, which conform to the conclusion in Corollary \ref{co-con}.
	\begin{figure}[h!]
		\centering
		{
			\includegraphics[width=0.52\textwidth]{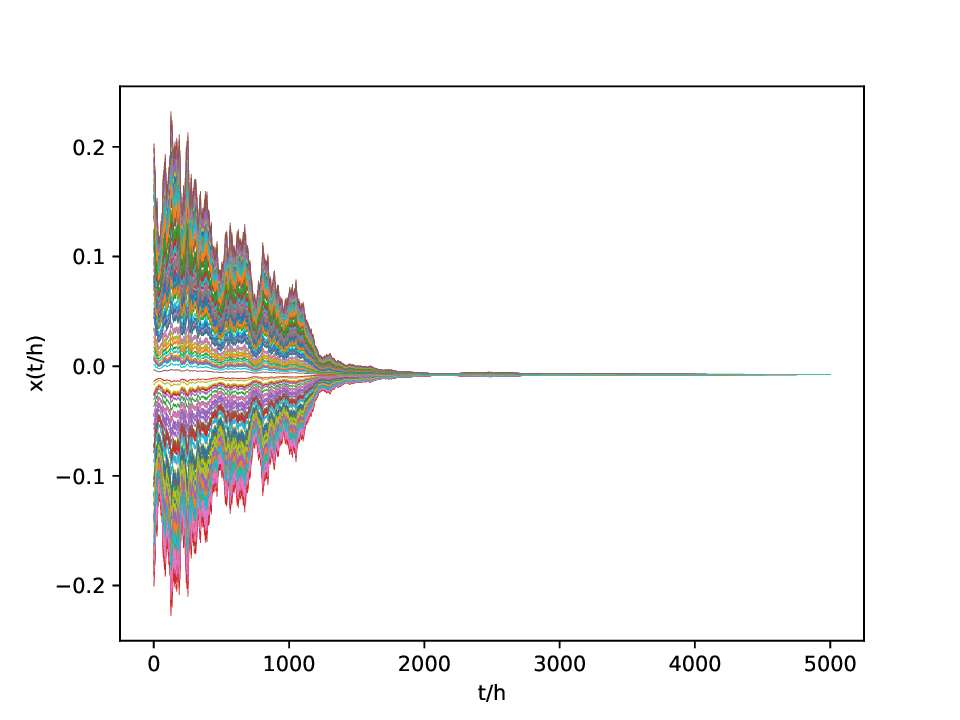}}
		\hspace{0in}
		\caption{The state trajectories of $150$ particles in Experiment 2 of Example \ref{example1}.}
		\label{fig3}
	\end{figure}
\end{example}

\begin{example}\label{example33}
	In this example, we illustrate the performance of the SCBO algorithm on some nonsmooth nonconvex test funtions. The chosen test functions to be $f$ in \eqref{problem} are provided in Table \ref{function}. Here, we choose the parameters randomly to show that the good performance of the proposed SCBO algorithm is not limited to the condition of parameters in Corollary \ref{co-con}. From the numerical results, we can see that the SCBO algorithm also performs well with high probability in this situation. 
	\begin{table}[h]
		\caption{Test functions in Example \ref{example33}.}
		\centering
		\renewcommand\arraystretch{2}
		\resizebox{\textwidth}{!}{
			\begin{tabular}{c|c} 
				\hline
				& Objective function $f$ in problem \eqref{problem} \\
				\hline
				$f_1$ & $\frac{1}{d}\sum_{l=1}^d[|x_l|-10\cos(2\pi x_l)+10]$ \\ 
				$f_2$ & $\sum_{l=1}^d|x_l|+\prod_{l=1}^d|x_l|$ \\ 
				$f_3$ & $\frac{1}{4000}\sum_{l=1}^d|x_l|-\prod_{l=1}^d\cos\left(\frac{x_l}{\sqrt{l}}\right)+1$ \\ 
				$f_4$ & $-10\exp\left( -0.2\sqrt{\frac{1}{d}\sum_{l=1}^d |x_l|}\right)-\exp\left(\frac{1}{d}\sum_{l=1}^d\cos(2\pi x_l) \right)+10+\exp(1)$ \\ 
				$f_5$ & $\left[\sum_{l=1}^d \sin^2(x_l)-\exp\left( -\sum_{l=1}^d x_l^2\right)\right]\exp\left( -\sum_{l=1}^d\sin^2(\sqrt{|x_l|})\right)+1$ \\ 
				\hline
		\end{tabular}}
		\label{function}
	\end{table}
	\begin{figure}[htbp]
		\centering
		\hspace{-0.37in}
		\subfigure[$f_1$]{
			\includegraphics[width=0.38\linewidth]{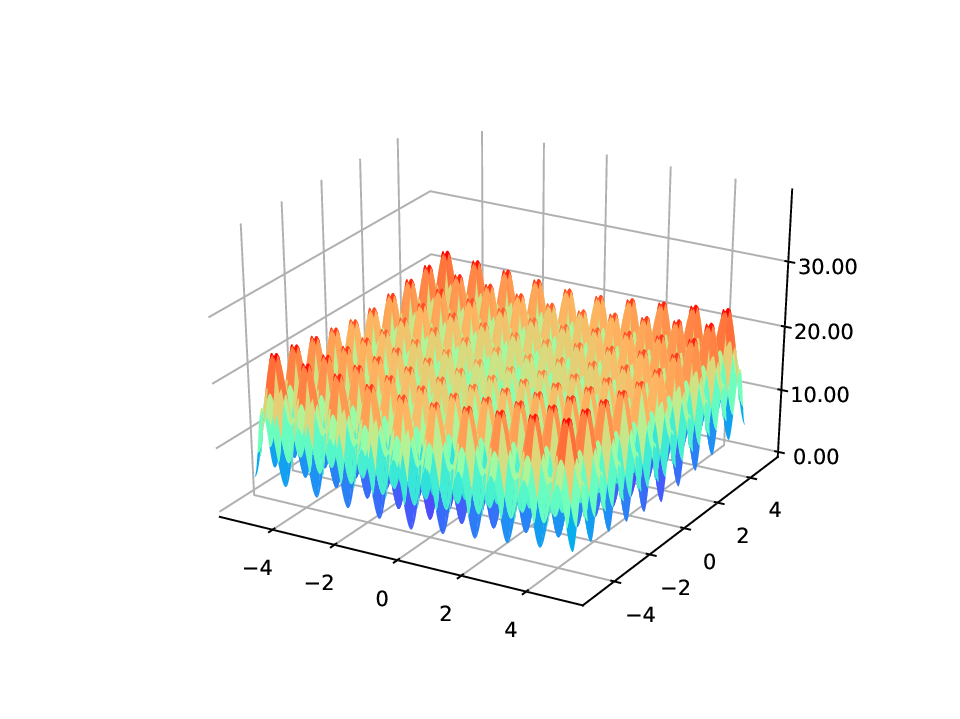}
		}
		\hspace{-0.37in}
		\subfigure[$f_2$]{
			\includegraphics[width=0.38\linewidth]{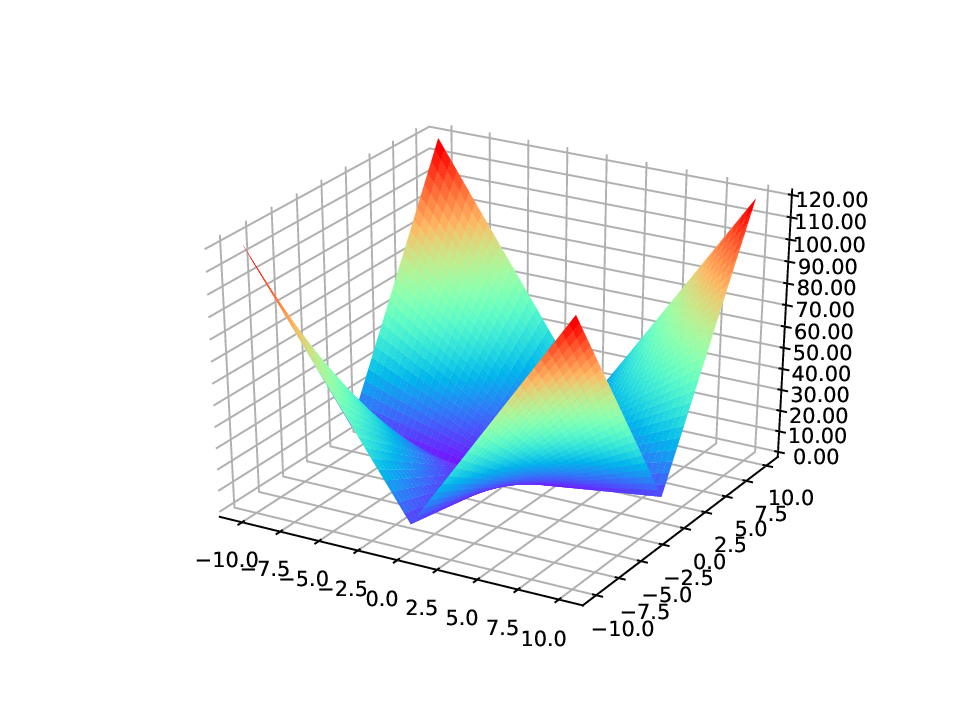}
		}
		\hspace{-0.326in}
		\subfigure[$f_3$]{
			\includegraphics[width=0.38\linewidth]{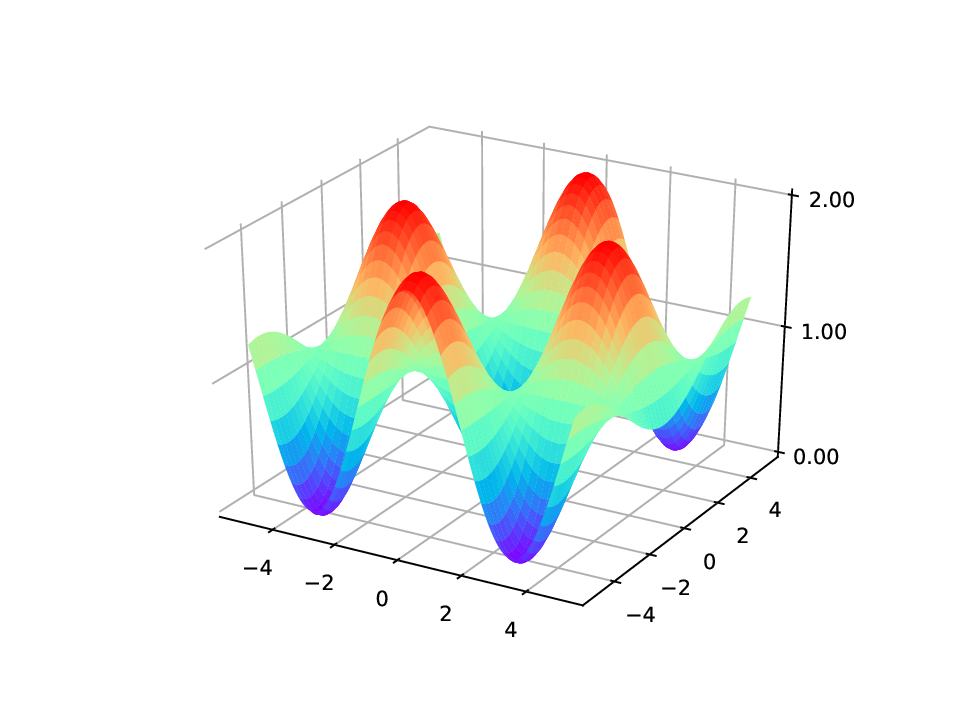}
		}
		\subfigure[$f_4$]{
			\includegraphics[width=0.38\linewidth]{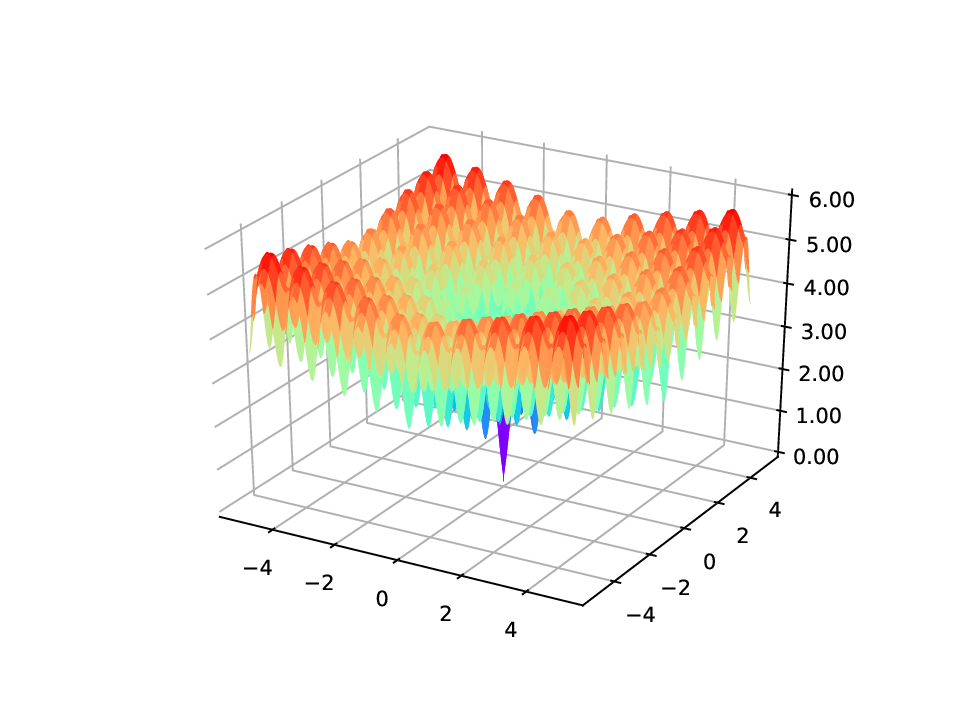}
		}
		\subfigure[$f_5$]{
			\includegraphics[width=0.38\linewidth]{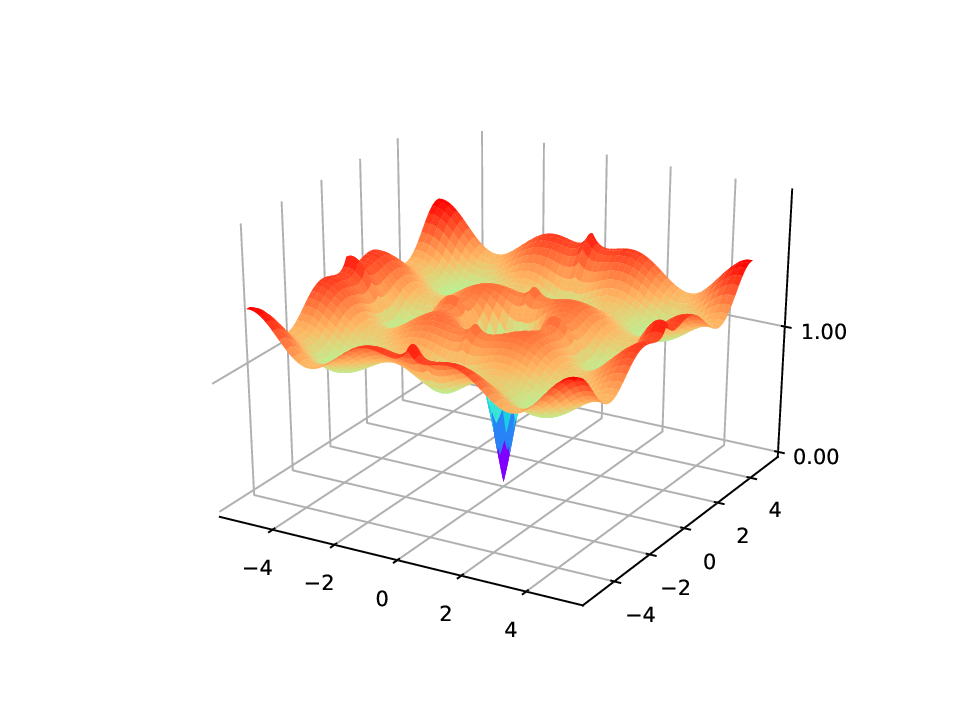}
		}
		\hspace{0.45in}		
		\caption{2D plots of the test functions in Table \ref{function} with $d=2$.} 
	\label{2Dfig}
\end{figure}

All these functions have a global minimum $f_{min}=0$ at $x^*=(0,\ldots,0)^{\mathrm{T}}\in\mathbb{R}^d$. We can use the function in  \eqref{smoothing} to construct a smoothing function with $q=0$ for each of these objective functions. Fig. \ref{2Dfig} provides the 2D plots of the five test functions in Table \ref{function} with $d=2$. In what follows, we always consider these functions with $d=2$.

\textbf{Experiment 1}: Show the unavailability of the deterministic smoothing algorithm, and explore the influence of parameters $N$ and $\beta$ in the SCBO algorithm.

We firstly apply the smoothing projected gradient (SPG) algorithm \cite{ZhangSmoothing2009} to solve the test function $f_1$. In the unconstrained setting, the SPG algorithm is formulated by
\begin{equation*}
	\begin{split}
		x_{k+1}&=x_k-\alpha_1 \nabla_x\tilde{f}_1(x_k,\mu_k),\\
		\mu_{k+1}&=\alpha_2\mu_k, \quad 0<\alpha_2<1,
	\end{split}
\end{equation*}
where a backtracking line search is performed to get the step length $\alpha_1$, and $\tilde{f}_1$ is the smoothing function of $f_1$ constructed by \eqref{smoothing}. 
We set the 100 uniformly distributed random points in $[-5,5]^2$ as the initial points. 
The numerical results of the SPG algorithm with these initial points are plotted in Fig. \ref{GD}, in which we refer the initial point that let the SPG algorithm find the global minima as a successful initial point, and the initial point that let the SPG get stuck at the local not global minima as a unsuccessful one. From Fig. \ref{GD}, we see that there is only one successful initial point, which is close enough to the optimal solution.  Thus, 
it can be seen that the SPG algorithm cannot solve this nonconvex test function efficiently.

\begin{figure}[h!]
	\centering
	{
		\includegraphics[width=0.52\textwidth]{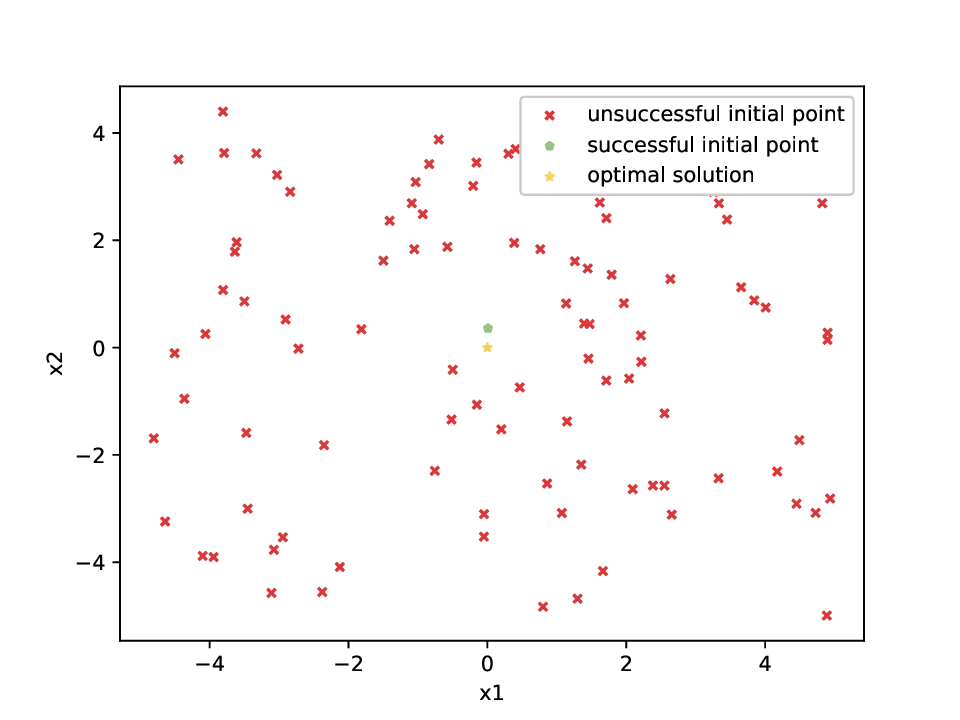}}
	\hspace{0in}
	\caption{Successful and unsuccessful initial points for SPG algorithm \cite{ZhangSmoothing2009}.}
	\label{GD}
\end{figure}

However, better results can be achieved when using the proposed CBO algorithm. We set $\mu_t=\mathrm{e}^{-0.9t}$, $\lambda=\sigma=1$ in the SCBO algorithm. The initial particles are also drawn from a uniform distribution within the domain $\mathcal{S}=[-5,5]^2$. We test this numerical simulation by $100$ realizations to get the success rates, and we call a realization successful if $f(x_\infty)$ is close to $f_{min}$, in the sense that
\begin{equation}\label{success}
	\frac{|f(x_\infty)-f_{min}|}{f_{max}-f_{min}}<0.005,
\end{equation}
where $f_{max}=\max_{x\in\mathcal{S}} f(x)$. Besides, the solution error of the SCBO algorithm is given by $\|x_\infty-x^*\|^2$ averaged over all runs.
Table \ref{tableN} records the success rates (rate) and the solution errors (sol-err) of the SCBO algorithm in terms of different particle numbers $N$ for different test functions with parameter $\beta=50$. And Table \ref{tablebeta} shows the results in terms of the values of $\beta$ for these test functions with particle number $N=40$. 
Fig. \ref{figureNbeta} presents the success rates in Table \ref{tableN} and Table \ref{tablebeta} visually. From these results, we find that more particles or larger $\beta$ can improve the performance of the SCBO algorithm on the whole. Except for functions $f_4$ and $f_5$ with small $N$ or $\beta$, the success rates for the other cases are good, while all the solution errors are promising. This is because the function values of both functions are relatively steep near the optimal solution, so even if the obtained point is very close to the optimal solution, the function value of this point cannot meet the condition for algorithm success. All these numerical results indicate that the SCBO algorithm is a promising global optimization method for solving nonsmooth nonconvex optimization problems.

\begin{table}[!h]
	\caption{Results of the SCBO algorithm with different particle numbers $N$.} \label{tableN}
	\begin{center}
		\renewcommand\arraystretch{1.2}
		\resizebox{\textwidth}{!}{
			\begin{tabular}{lllllllll} \hline
				&$N$&$20$&$50$&$80$&$110$&$140$&$170$&$200$\\
				\hline
				\multirow{2}{*}{$f_1$}& rate&$91\%$&$95\%$&$96\%$&$100\%$&$100\%$&$100\%$&$100\%$\\
				&sol-err &$1.44\mathrm{e}$-$05$&$1.55\mathrm{e}$-$05$&$7.32\mathrm{e}$-$06$&$2.60\mathrm{e}$-$06$&$2.12\mathrm{e}$-$06$&$8.38\mathrm{e}$-$07$&$4.26\mathrm{e}$-$07$\\
				\hline
				\multirow{2}{*}{$f_2$}&rate&$97\%$&$97\%$&$100\%$&$100\%$&$100\%$&$100\%$&$100\%$\\
				&sol-err&$3.04\mathrm{e}$-$04$&$4.58\mathrm{e}$-$05$&$2.06\mathrm{e}$-$05$&$1.08\mathrm{e}$-$05$&$6.93\mathrm{e}$-$06$&$6.54\mathrm{e}$-$06$&$9.42\mathrm{e}$-$07$\\
				\hline
				\multirow{2}{*}{$f_3$}&rate&$79\%$&$84\%$&$81\%$&$96\%$&$96\%$&$97\%$&$100\%$\\
				&sol-err&$8.26\mathrm{e}$-$03$&$7.45\mathrm{e}$-$03$&$6.42\mathrm{e}$-$03$&$2.31\mathrm{e}$-$03$&$2.54\mathrm{e}$-$03$&$1.25\mathrm{e}$-$03$&$1.04\mathrm{e}$-$03$\\
				\hline
				\multirow{2}{*}{$f_4$}&rate&$41\%$&$53\%$&$52\%$&$58\%$&$79\%$&$88\%$&$95\%$\\
				&sol-err&$1.12\mathrm{e}$-$06$&$7.56\mathrm{e}$-$07$&$4.43\mathrm{e}$-$07$&$4.55\mathrm{e}$-$07$&$1.88\mathrm{e}$-$07$&$6.31\mathrm{e}$-$08$&$3.02\mathrm{e}$-$08$\\
				\hline
				\multirow{2}{*}{$f_5$}&rate&$44\%$&$39\%$&$55\%$&$65\%$&$73\%$&$75\%$&$97\%$\\
				&sol-err&$9.79\mathrm{e}$-$05$&$2.89\mathrm{e}$-$04$&$9.73\mathrm{e}$-$05$&$3.44\mathrm{e}$-$05$&$1.44\mathrm{e}$-$05$&$1.34\mathrm{e}$-$05$&$1.12\mathrm{e}$-$05$\\
				\hline
		\end{tabular}}
	\end{center}
\end{table}

\begin{table}[!h]
	\caption{Results of the SCBO algorithm with different values of $\beta$.} \label{tablebeta}
	\begin{center}
		\renewcommand\arraystretch{1.2}
		\resizebox{\textwidth}{!}{
			\begin{tabular}{lllllllll} \hline
				&$\beta$&$20$&$40$&$60$&$80$&$100$&$120$&$140$\\
				\hline
				\multirow{2}{*}{$f_1$}&rate&$78\%$&$83\%$&$90\%$&$98\%$&$97\%$&$100\%$&$100\%$\\
				&sol-err&$5.41\mathrm{e}$-$04$&$4.59\mathrm{e}$-$05$&$1.61\mathrm{e}$-$05$&$8.74\mathrm{e}$-$06$&$6.39\mathrm{e}$-$07$&$2.65\mathrm{e}$-$06$&$4.22\mathrm{e}$-$07$\\
				\hline
				\multirow{2}{*}{$f_2$}&rate&$92\%$&$97\%$&$98\%$&$100\%$&$100\%$&$100\%$&$100\%$\\
				&sol-err&$1.07\mathrm{e}$-$03$&$1.01\mathrm{e}$-$03$&$1.78\mathrm{e}$-$04$&$2.26\mathrm{e}$-$05$&$3.29\mathrm{e}$-$05$&$1.99\mathrm{e}$-$05$&$4.96\mathrm{e}$-$06$\\
				\hline
				\multirow{2}{*}{$f_3$}&rate&$85\%$&$90\%$&$89\%$&$100\%$&$100\%$&$100\%$&$100\%$\\
				&sol-err&$6.17\mathrm{e}$-$03$&$3.36\mathrm{e}$-$03$&$7.51\mathrm{e}$-$03$&$4.75\mathrm{e}$-$03$&$3.64\mathrm{e}$-$03$&$2.97\mathrm{e}$-$03$&$1.13\mathrm{e}$-$03$\\
				\hline
				\multirow{2}{*}{$f_4$}&rate&$21\%$&$33\%$&$51\%$&$47\%$&$76\%$&$96\%$&$98\%$\\
				&sol-err&$2.63\mathrm{e}$-$06$&$5.11\mathrm{e}$-$07$&$1.50\mathrm{e}$-$07$&$3.40\mathrm{e}$-$07$&$1.34\mathrm{e}$-$08$&$5.43\mathrm{e}$-$09$&$6.11\mathrm{e}$-$09$\\
				\hline
				\multirow{2}{*}{$f_5$}&rate&$19\%$&$24\%$&$30\%$&$39\%$&$42\%$&$93\%$&$93\%$\\
				&sol-err&$1.05\mathrm{e}$-$04$&$3.76\mathrm{e}$-$05$&$2.92\mathrm{e}$-$05$&$1.65\mathrm{e}$-$05$&$2.00\mathrm{e}$-$05$&$7.83\mathrm{e}$-$06$&$5.86\mathrm{e}$-$06$\\
				\hline
		\end{tabular}}
	\end{center}
\end{table}

\begin{figure}[!]
	\centering
	\subfigure[]{
		\includegraphics[width=0.5\linewidth]{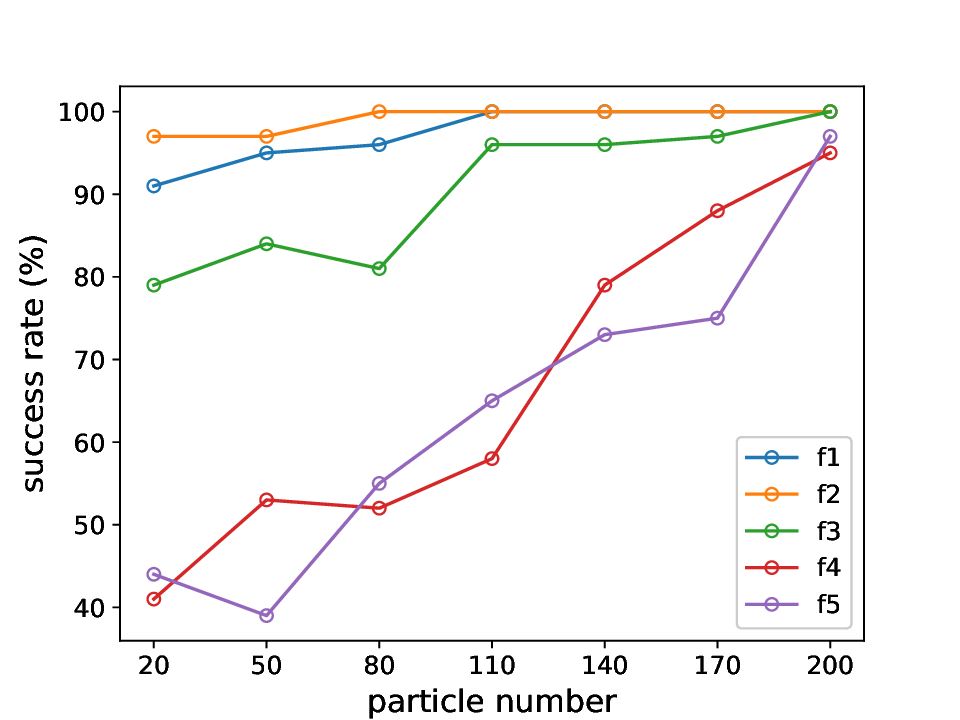}
	}
	\hspace{-0.3in}
	\subfigure[]{
		\includegraphics[width=0.5\linewidth]{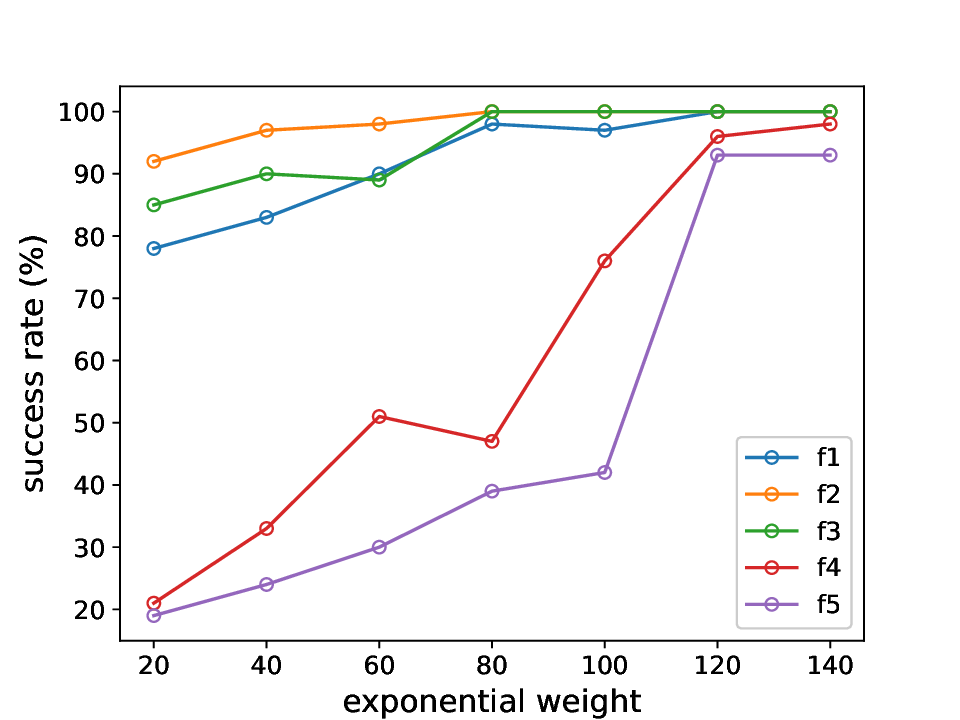}
	}
	\caption{$(\mathrm{a})$: Variations of success rates with different particle number $N$; $(\mathrm{b})$: Variations of success rates with different exponential weight $\beta$.}
	\label{figureNbeta}
\end{figure}

	\textbf{Experiment 2}: Comparison with the CBO algorithm in \cite{HaConvergence2020}.
	
	We compare the performance of the proposed SCBO algorithm with the CBO algorithm in \cite{HaConvergence2020} by considering the test functions in Table \ref{function}. For the both algorithms, we set the same initial datas, which are the uniformly distributed random variables in $\mathcal{S}=[-5,5]^d$, and the same parameters $\lambda=1$, $\sigma=0.5$. Although the SCBO algorithm is related to the constructed smoothing function of $f$, its numerical performance is not limited by the choice of smoothing function. Therefore, in this experiment, we attempt to choose the function in \eqref{smoothing2}
	as the smoothing function of nonsmooth term $|x_l|$ to construct a smoothing function of these nonsmooth functions. Set the update of the smoothing parameter by $\mu_t=0.05\mathrm{e}^{-0.9t}$. We still adopt \eqref{success} as the criterion for the success of this experiment. Tables \ref{com1} and \ref{com2} report the success rates (rate), objective function values (fun-val) and solution errors (sol-err) given by $\|x_\infty-x^*\|^2 $, averaged over 100 runs. Specifically, Table \ref{com1} shows the comparisons between the two methods with a fixed $\beta=40$ and different values of $N\in\{100,150,200\}$, and Table \ref{com2} shows the comparisons with a fixed $N=90$ and different values of $\beta\in\{80,100,120\}$. The results illustrate that the performance of SCBO algorithm is on par with the CBO algorithm \cite{HaConvergence2020}, and even slightly better than the CBO algorithm in some cases. From the theoretical results of this paper and the comparisons on numerical results, we can conclude that the SCBO algorithm owns some theoretical improvements over the CBO algorithm, while ensuring on par or better numerical performance. 
\begin{table}[!h]
	\caption{Comparisons between SCBO and CBO in \cite{HaConvergence2020} with different values of $N$.} \label{com1}
	\begin{center}
		\renewcommand\arraystretch{1.2}
		\resizebox{\textwidth}{!}{
			\begin{tabular}{llllllll} \hline
				&&\multicolumn {3}{c}{CBO}&\multicolumn {3}{c}{SCBO}\\
				\cmidrule(lr{0pt}){3-5}\cmidrule(lr{0pt}){6-8}
				&&$N=100$&$N=150$&$N=200$&$N=100$&$N=150$&$N=200$\\
				\hline
				$f_1$&rate&$\textbf{100\%}$&$\textbf{100\%}$&$\textbf{100\%}$&$\textbf{100\%}$&$98\%$&$\textbf{100\%}$\\
				~&fun-val&$4.79\mathrm{e}$-$03$
				&$\textbf{3.17}\mathrm{\textbf{e}}$\bf{-}$\textbf{03}$&$\textbf{3.02}\mathrm{\textbf{e}}$\bf{-}$\textbf{03}$&$\textbf{3.96}\mathrm{\textbf{e}}$\bf{-}$\textbf{03}$&$7.84\mathrm{e}$-$03$&$3.09\mathrm{e}$-$03$\\
				~&sol-err&$2.41\mathrm{e}$-$05$
				&$\textbf{1.37}\mathrm{\textbf{e}}$\bf{-}$\textbf{05}$&$1.58\mathrm{e}$-$05$&$\textbf{1.76}\mathrm{\textbf{e}}$\bf{-}$\textbf{05}$&$5.31\mathrm{e}$-$05$&$\textbf{1.21}\mathrm{\textbf{e}}$\bf{-}$\textbf{05}$\\
				\hline
				$f_2$&rate&$\textbf{100\%}$&$\textbf{100\%}$&$\textbf{100\%}$&$\textbf{100\%}$&$\textbf{100\%}$&$\textbf{100\%}$\\
				~&fun-val&$\textbf{5.80}\mathrm{\textbf{e}}$\bf{-}$\textbf{03}$
				&$\textbf{4.72}\mathrm{\textbf{e}}$\bf{-}$\textbf{03}$&$4.16\mathrm{e}$-$03$&$7.17\mathrm{e}$-$03$&$5.14\mathrm{e}$-$03$&$\textbf{4.05}\mathrm{\textbf{e}}$\bf{-}$\textbf{03}$\\
				~&sol-err&$\textbf{3.16}\mathrm{\textbf{e}}$\bf{-}$\textbf{05}$
				&$\textbf{1.92}\mathrm{\textbf{e}}$\bf{-}$\textbf{05}$&$\textbf{1.80}\mathrm{\textbf{e}}$\bf{-}$\textbf{05}$&$6.75\mathrm{e}$-$05$&$4.95\mathrm{e}$-$05$&$2.45\mathrm{e}$-$05$\\
				\hline
				$f_3$&rate&$\textbf{92\%}$&$\textbf{94\%}$&$94\%$&$89\%$&$92\%$&$\textbf{96\%}$\\
				~&fun-val&$\textbf{6.57}\mathrm{\textbf{e}}$\textbf{-}$\textbf{03}$
				&$4.94\mathrm{e}$-$03$&$2.54\mathrm{e}$-$03$&$7.08\mathrm{e}$-$03$&$\textbf{4.33}\mathrm{\textbf{e}}$\textbf{-}$\textbf{03}$&$\textbf{1.84}\mathrm{\textbf{e}}$\textbf{-}$\textbf{03}$\\
				~&sol-err&$3.11\mathrm{e}$-$01$
				&$\textbf{2.45}\mathrm{\textbf{e}}$\textbf{-}$\textbf{01}$&$8.54\mathrm{e}$-$03$&$\textbf{2.75}\mathrm{\textbf{e}}$\textbf{-}$\textbf{01}$&$2.57\mathrm{e}$-$01$&$\textbf{6.06}\mathrm{\textbf{e}}$\textbf{-}$\textbf{03}$\\
				\hline
				$f_4$&rate&$54\%$&$\textbf{63\%}$&$\textbf{69\%}$&$\textbf{66\%}$&$60\%$&$67\%$\\
				~&fun-val&$4.67\mathrm{e}$-$02$
				&$\textbf{3.56}\mathrm{\textbf{e}}$\textbf{-}$\textbf{02}$&$3.17\mathrm{e}$-$02$&$\textbf{3.51}\mathrm{\textbf{e}}$\textbf{-}$\textbf{02}$&$3.81\mathrm{e}$-$02$&$\textbf{3.08}\mathrm{\textbf{e}}$\textbf{-}$\textbf{02}$\\
				~&sol-err&$5.38\mathrm{e}$-$05$
				&$\textbf{3.03}\mathrm{\textbf{e}}$\textbf{-}$\textbf{06}$&$1.23\mathrm{e}$-$05$&$\textbf{9.88}\mathrm{\textbf{e}}$\textbf{-}$\textbf{07}$&$4.54\mathrm{e}$-$06$&$\textbf{5.00}\mathrm{\textbf{e}}$\textbf{-}$\textbf{07}$\\
				\hline
				$f_5$&rate&$67\%$&$\textbf{82\%}$&$\textbf{89\%}$&$\textbf{72\%}$&$\textbf{82\%}$&$86\%$\\
				~&fun-val&$7.07\mathrm{e}$-$03$
				&$4.95\mathrm{e}$-$03$&$4.67\mathrm{e}$-$03$&$\textbf{6.48}\mathrm{\textbf{e}}$\textbf{-}$\textbf{03}$&$\textbf{4.91}\mathrm{\textbf{e}}$\textbf{-}$\textbf{03}$&$\textbf{4.20}\mathrm{\textbf{e}}$\textbf{-}$\textbf{03}$\\
				~&sol-err&$7.97\mathrm{e}$-$05$
				&$\textbf{2.52}\mathrm{\textbf{e}}$\textbf{-}$\textbf{05}$&$2.41\mathrm{e}$-$05$&$\textbf{3.94}\mathrm{\textbf{e}}$\textbf{-}$\textbf{05}$&$2.79\mathrm{e}$-$05$&$\textbf{1.55}\mathrm{\textbf{e}}$\textbf{-}$\textbf{05}$\\
				\hline
		\end{tabular}}
\end{center}
\end{table}

\begin{table}[!h]
	\caption{Comparisons between SCBO and CBO in \cite{HaConvergence2020} with different values of $\beta$.} \label{com2}
	\begin{center}
		\renewcommand\arraystretch{1.2}
		\resizebox{\textwidth}{!}{
			\begin{tabular}{llllllll} \hline
				&&\multicolumn {3}{c}{CBO}&\multicolumn {3}{c}{SCBO}\\
				\cmidrule(lr{0pt}){3-5}\cmidrule(lr{0pt}){6-8}
				&&$\beta=80$&$\beta=100$&$\beta=120$&$\beta=80$&$\beta=100$&$\beta=120$\\
				\hline
				$f_1$&rate&$\textbf{98\%}$&$\textbf{100\%}$&$\textbf{100\%}$&$96\%$&$99\%$&$\textbf{100\%}$\\
				~&fun-val&$\textbf{2.64}\mathrm{\textbf{e}}$\textbf{-}$\textbf{02}$
				&$\textbf{5.38}\mathrm{\textbf{e}}$\textbf{-}$\textbf{03}$&$\textbf{3.92}\mathrm{\textbf{e}}$\textbf{-}$\textbf{03}$&$3.98\mathrm{e}$-$02$&$7.35\mathrm{e}$-$03$&$4.86\mathrm{e}$-$03$\\
				~&sol-err&$\textbf{1.99}\mathrm{\textbf{e}}$\textbf{-}$\textbf{02}$
				&$\textbf{3.94}\mathrm{\textbf{e}}$\textbf{-}$\textbf{05}$&$4.32\mathrm{e}$-$05$&$2.08\mathrm{e}$-$02$&$6.98\mathrm{e}$-$05$&$\textbf{2.81}\mathrm{\textbf{e}}$\textbf{-}$\textbf{05}$\\
				\hline
				$f_2$&rate&$\textbf{100\%}$&$\textbf{100\%}$&$\textbf{100\%}$&$\textbf{100\%}$&$\textbf{100\%}$&$\textbf{100\%}$\\
				~&fun-val&$\textbf{4.43}\mathrm{\textbf{e}}$\textbf{-}$\textbf{03}$
				&$4.20\mathrm{e}$-$03$&$3.87\mathrm{e}$-$03$&$5.07\mathrm{e}$-$03$&$\textbf{4.08}\mathrm{\textbf{e}}$\textbf{-}$\textbf{03}$&$\textbf{3.80}\mathrm{\textbf{e}}$\textbf{-}$\textbf{03}$\\
				~&sol-err&$\textbf{1.89}\mathrm{\textbf{e}}$\textbf{-}$\textbf{05}$
				&$\textbf{2.35}\mathrm{\textbf{e}}$\textbf{-}$\textbf{05}$&$4.31\mathrm{e}$-$05$&$5.53\mathrm{e}$-$05$&$5.96\mathrm{e}$-$05$&$\textbf{2.99}\mathrm{\textbf{e}}$\textbf{-}$\textbf{05}$\\
				\hline
				$f_3$&rate&$89\%$&$\textbf{96\%}$&$90\%$&$\textbf{92\%}$&$95\%$&$\textbf{92\%}$\\
				~&fun-val&$5.94\mathrm{e}$-$03$
				&$\textbf{3.22}\mathrm{\textbf{e}}$\textbf{-}$\textbf{03}$&$\textbf{3.39}\mathrm{\textbf{e}}$\textbf{-}$\textbf{03}$&$\textbf{3.33}\mathrm{\textbf{e}}$\textbf{-}$\textbf{03}$&$4.19\mathrm{e}$-$03$&$6.07\mathrm{e}$-$03$\\
				~&sol-err&$1.98\mathrm{e}$-$02$
				&$6.79\mathrm{e}$-$01$&$\textbf{7.03}\mathrm{\textbf{e}}$\textbf{-}$\textbf{01}$&$\textbf{1.14}\mathrm{\textbf{e}}$\textbf{-}$\textbf{02}$&$\textbf{6.50}\mathrm{\textbf{e}}$\textbf{-}$\textbf{01}$&$7.29\mathrm{e}$-$01$\\
				\hline
				$f_4$&rate&$55\%$&$\textbf{71\%}$&$72\%$&$\textbf{58\%}$&$65\%$&$\textbf{76\%}$\\
				~&fun-val&$4.89\mathrm{e}$-$02$
				&$\textbf{3.38}\mathrm{\textbf{e}}$\textbf{-}$\textbf{02}$&$3.22\mathrm{e}$-$02$&$\textbf{3.99}\mathrm{\textbf{e}}$\textbf{-}$\textbf{02}$&$3.72\mathrm{e}$-$02$&$\textbf{3.06}\mathrm{\textbf{e}}$\textbf{-}$\textbf{02}$\\
				~&sol-err&$1.56\mathrm{e}$-$04$
				&$5.38\mathrm{e}$-$06$&$1.26\mathrm{e}$-$05$&$\textbf{1.36}\mathrm{\textbf{e}}$\textbf{-}$\textbf{04}$&$\textbf{4.46}\mathrm{\textbf{e}}$\textbf{-}$\textbf{06}$&$\textbf{1.01}\mathrm{\textbf{e}}$\textbf{-}$\textbf{05}$\\
				\hline
				$f_5$&rate&$\textbf{90\%}$&$88\%$&$93\%$&$88\%$&$\textbf{95\%}$&$\textbf{97\%}$\\
				~&fun-val&$5.13\mathrm{e}$-$03$
				&$5.03\mathrm{e}$-$03$&$\textbf{3.23}\mathrm{\textbf{e}}$\textbf{-}$\textbf{03}$&$\textbf{4.51}\mathrm{\textbf{e}}$\textbf{-}$\textbf{03}$&$\textbf{3.75}\mathrm{\textbf{e}}$\textbf{-}$\textbf{03}$&$3.44\mathrm{e}$-$03$\\
				~&sol-err&$5.20\mathrm{e}$-$05$
				&$7.82\mathrm{e}$-$05$&$\textbf{1.30}\mathrm{\textbf{e}}$\textbf{-}$\textbf{05}$&$\textbf{3.59}\mathrm{\textbf{e}}$\textbf{-}$\textbf{05}$&$\textbf{3.28}\mathrm{\textbf{e}}$\textbf{-}$\textbf{05}$&$2.87\mathrm{e}$-$05$\\
				\hline
		\end{tabular}}
\end{center}
\end{table}
\end{example}

\section{Conclusions}

\label{sec:6}
In this paper, we consider a unconstrained optimization problem with a nonsmooth nonconvex objective function and the convergence analysis of the finite particle system of a variant of the CBO method.
First, we propose a variant of CBO algorithm based on a proper smoothing function.
Second, using Ito's formula, we derive an explicit formula of the state differences of different particles.
This exact formula yields that the proposed algorithm exhibits a global consensus almost surely.
Moreover, it is proved that under certain conditions on the drift and diffusion parameters, there exists a common consensus state which is the almost-sure limit of the solution of the proposed algorithm.
Finally, inspired by the smoothing techniques, we give a sufficient condition on the system parameters and initial data to
guarantee that the function value at global consensus state can be close enough to the global minimum.
And note that the conditions on the parameters are independent of the dimension of the problem.
Moreover, besides the theoretical improvement, this work also provides a good starting point for studying the convergence of the original CBO algorithm under weaker assumptions and other metaheuristic algorithms based on the swarm intelligence. Future studies could fruitfully explore these issues.

\begin{acknowledgements}
The research of this author is partially supported by National Key Research and Development Program of China (2021YFA1003500), National Natural Science Foundation of China grants (12271127, 62176073) and Fundamental Research Funds for Central Universities (HIT.OCEF.2024050, 2022FRFK060017). 
\end{acknowledgements}

\section*{Data Availability}
The datasets generated during the current study are available from the corresponding author on reasonable request.
%
%
%
%
%
%
%
%
%
\bibliographystyle{spmpsci} 
\bibliography{cbosmoothing}

\begin{thebibliography}{10}
\providecommand{\url}[1]{{#1}}
\providecommand{\urlprefix}{URL }
\expandafter\ifx\csname urlstyle\endcsname\relax
  \providecommand{\doi}[1]{DOI~\discretionary{}{}{}#1}\else
  \providecommand{\doi}{DOI~\discretionary{}{}{}\begingroup
  \urlstyle{rm}\Url}\fi

\bibitem{BerntStochastic1985}
Bernt, O.: Stochastic Differential Equations: An Introduction with
  Applications.
\newblock Springer, Berlin, Heridelberg (1985)

\bibitem{BianWorst2013}
Bian, W., Chen, X.J.: Worst-case complexity of smoothing quadratic
  regularization methods for non-{L}ipschitzian optimization.
\newblock SIAM J. Optim. \textbf{23}(3), 1718--1741 (2013)

\bibitem{BianNeural2014}
Bian, W., Chen, X.J.: Neural network for nonsmooth, nonconvex constrained
  minimization via smooth approximation.
\newblock IEEE Trans. Neural Netw. Learn. Syst. \textbf{25}(3), 545--556 (2014)

\bibitem{BianA2020}
Bian, W., Chen, X.J.: A smoothing proximal gradient algorithm for nonsmooth
  convex regression with cardinality penalty.
\newblock SIAM J. Numer. Anal. \textbf{58}(1), 858--883 (2020)

\bibitem{BorghiConstrained2023}
Borghi, G., Herty, M., Pareschi, L.: Constrained consensus-based optimization.
\newblock SIAM J. Optim. \textbf{33}(1), 211--236 (2023)

\bibitem{CarrilloAn2016}
Carrillo, J.A., Choi, Y.P., Totzeck, C., Tse, O.: An analytical framework for a
  consensus-based global optimization method.
\newblock Math. Models Meth. Appl. Sci. \textbf{28}, 1037--1066 (2018)

\bibitem{CarrilloA2019}
Carrillo, J.A., Jin, S., Li, L., Zhu, Y.: A consensus-based global optimization
  method for high dimensional machine learning problems.
\newblock ESAIM Control Optim. Calc. Var. \textbf{27}(S5) (2021)

\bibitem{JoseFedCBO2023}
Carrillo, J.A., Trillos, N.G., Li, S., Zhu, Y.: Fed{CBO}: {R}eaching group
  consensus in clustered federated learning through consensus-based
  optimization.
\newblock arXiv preprint arXiv: 2305.02894 (2023)

\bibitem{ChenSmoothing2012}
Chen, X.J.: Smoothing methods for nonsmooth, nonconvex minimization.
\newblock Math. Program. \textbf{134}(1), 71--99 (2012)

\bibitem{ChenA2018}
Chen, X.J., Kelley, C.T., Xu, F.M., Zhang, Z.K.: A smoothing direct search
  method for {M}onte {C}arlo-based bound constrained composite nonsmooth
  optimization.
\newblock SIAM J. Sci. Comput. \textbf{40}(4), A2174--A2199 (2018)

\bibitem{ChenPenalty2016}
Chen, X.J., Lu, Z.S., Pong, T.K.: Penalty methods for a class of
  non-{L}ipschitz optimization problems.
\newblock SIAM J. Optim. \textbf{26}(3), 1465--1492 (2016)

\bibitem{Chenchene2023}
Chenchene, E., Huang, H., Qiu, J.: A consensus-based algorithm for non-convex
  multiplayer games.
\newblock arXiv preprint arXiv: 2311.08270 (2023)

\bibitem{CrowLognormal1988}
Crow, E.L., Shimizu, K.: Lognormal Distributions: Theory and Application.
\newblock Marcel-Dekker Inc., New York (1988)

\bibitem{DemboLarge1998}
Dembo, A., Zeitouni, O.: Large Deviations Techniques and Applications.
\newblock Springer, Berlin, Heidelberg (1998)

\bibitem{DemoA2021}
Demo, N., Tezzele, M., Rozza, G.: A supervised learning approach involving
  active subspaces for an efficient genetic algorithm in high-dimensional
  optimization problems.
\newblock SIAM J. Sci. Comput. \textbf{43}(3), B831--B853 (2021)

\bibitem{FornasierConsenseus2021}
Fornasier, M., Huang, H., Pareschi, L., Sunnen, P.: Consensus-based
  optimization on the sphere {II}: Convergence to global minimizers and machine
  learning.
\newblock J. Mach. Learn. Res. \textbf{22}(237), 1--55 (2021)

\bibitem{FornasierConsensus2022}
Fornasier, M., Klock, T., Riedl, K.: Consensus-based optimization methods
  converge globally.
\newblock arXiv preprint arXiv: 2103.15130 (2022)

\bibitem{FornasierConsensus2024}
Fornasier, M., Richtarik, P., Ried, K., Sun, L.: Consensus-based optimization
  with truncated noise.
\newblock arXiv preprint arXiv: 2310.16610 (2024)

\bibitem{2013Brownian}
Gall, J.: Brownian Motion, Martingales, and Stochastic Calculus.
\newblock Springer, Berlin (2013)

\bibitem{Grefenstette1986Optimization}
Grefenstette, J.J.: Optimization of control parameters for genetic algorithms.
\newblock IEEE Trans. Syst. Man Cybern. -Syst. \textbf{16}(1), 122--128 (1986)

\bibitem{HaConvergence2020}
Ha, S.Y., Jin, S., Kim, D.: Convergence of a first-order consensus-based global
  optimization algorithm.
\newblock Math. Models Meth. Appl. Sci. \textbf{30}, 2417--2444 (2020)

\bibitem{HaConvergence2021}
Ha, S.Y., Jin, S., Kim, D.: Convergence and error estimates for time-discrete
  consensus-based optimization algorithms.
\newblock Numer. Math. \textbf{147}, 255--282 (2021)

\bibitem{HertleinAn2019}
Hertlein, L., Ulbrich, M.: An inexact bundle algorithm for nonconvex nonsmooth
  minimization in {H}ilbert space.
\newblock SIAM J. Control Optim. \textbf{57}(5), 3137--3165 (2019)

\bibitem{HiriartConvex1993}
Hiriart-Urruty, J.B., Lemarechal, C.: Convex Analysis and Minimization
  Algorithms I.
\newblock Springer-Verlag, Berlin (1993)

\bibitem{Huang2023}
Huang, H., Qiu, J., Riedl, K.: Consensus-based optimization for saddle point
  problems.
\newblock arXiv preprint arXiv: 2212.12334 (2023)

\bibitem{HuangRobust2018}
Huang, J., Jiao, Y.L., Lu, X.L., Zhu, L.P.: Robust decoding from 1-bit
  compressive sampling with ordinary and regularized least squares.
\newblock SIAM J. Sci. Comput. \textbf{40}(4), A2062--A2086 (2018)

\bibitem{Hwang1988Simulated}
Hwang, C.R.: Simulated annealing: Theory and applications.
\newblock Acta Appl. Math. \textbf{12}(1), 108--111 (1988)

\bibitem{AchimProbability2008}
Klenke, A.: Probability Theory.
\newblock Springer, London (2008)

\bibitem{MekaRank2008}
Meka, R., Jain, P., Caramanis, C., Dhillon, I.S.: Rank minimization via online
  learning.
\newblock In: Proceedings of the 25th International Conference on Machine
  Learning. Helsinki, Finland (2008)

\bibitem{GrigoriosStochastic2014}
Pavliotis, G.A.: Stochastic Processes and Applications.
\newblock Springer, NewYork (2014)

\bibitem{PinnauA2016}
Pinnau, R., Totzeck, C., Tse, O., Martin, S.: A consensus-based model for
  global optimization and its mean-field limit.
\newblock Math. Models Meth. Appl. Sci. \textbf{27}(1), 183--204 (2016)

\bibitem{RockafellarVariational1998}
Rockafellar, R., Wets, R.B.: Variational Analysis.
\newblock Springer, New York (1998)

\bibitem{RossStochastic1995}
Ross, S.M.: Stochastic Processes.
\newblock John Wiley \& Sons, New York (1995)

\bibitem{RuszA2021}
Ruszczyński, A.: A stochastic subgradient method for nonsmooth nonconvex
  multilevel composition optimization.
\newblock SIAM J. Control Optim. \textbf{59}(3), 2301--2320 (2021)

\bibitem{SiaryTabu2000}
Siarry, C.P.: Tabu search applied to global optimization.
\newblock European J. Oper. Res. \textbf{123}(2), 256--270 (2000)

\bibitem{Tarasewich2002Swarm}
Tarasewich, P., Mcmullen, P.R.: Swarm intelligence.
\newblock Commun. ACM \textbf{45}(8), 62--67 (2002)

\bibitem{2020Consensus}
Totzeck, C., Wolfram, M.T.: Consensus-based global optimization with personal
  best.
\newblock Math. Biosci. Eng. \textbf{17}, 6026--6044 (2020)

\bibitem{YangSwarm2018}
Yang, X.S., Deb, S., Zhao, Y.X., Fong, S., He, X.S.: Swarm intelligence:
  {P}ast, present and future.
\newblock Soft Comput. \textbf{22}, 5923--5933 (2018)

\bibitem{ZhangSmoothing2009}
Zhang, C., Chen, X.J.: Smoothing projected gradient method and its application
  to stochastic linear complementarity problems.
\newblock SIAM J. Optim. \textbf{20}(2), 627--649 (2009)

\end{thebibliography}
\end{document}